\documentclass[twocolumn,aps,pre,floatfix,superscriptaddress,longbibliography]{revtex4-1}
\pdfoutput=1 
\usepackage{amsmath,amssymb,eucal,graphicx,float,epstopdf,xparse}
\usepackage{epsfig,subfigure}
\usepackage[utf8]{inputenc}

\NewDocumentCommand{\eulerian}{omm}
 {%
  \genfrac<>{0pt}{}{#2}{#3}%
  \IfValueT{#1}{_{\!#1}}%
 }

\usepackage[colorlinks=true, urlcolor=blue, anchorcolor=blue, citecolor=blue,filecolor=blue,linkcolor=blue,menucolor=blue]{hyperref}

\begin{document}
\title{Random Maps with Sociological Flavor} 

\author{P. L. Krapivsky}
\affiliation{Department of Physics, Boston University, Boston, Massachusetts 02215, USA}
\affiliation{Santa Fe Institute, Santa Fe, New Mexico 87501, USA}

\begin{abstract} 
A map of a set to itself admits a representation by a graph with vertices being the elements of the set and an edge between every vertex and its image. Communities defined as the maximal connected components are uni-cyclic. The distributions of the sizes of communities and lengths of cycles for unconstrained random maps is a classical subject. We call experts the images and followers the remaining vertexes, and we further define prophets, egocentrics, and introverts. We introduce and analyze classes of random maps with sociological flavor. 
\end{abstract}

\maketitle

\section{Introduction}
\label{sec:intro}

A map of a set $S$ to itself, $f: S\to S$, can be represented by a directed graph in which elements of $S$ are vertices and a directed edge goes from each $x\in S$  to $f(x)$. This graph decomposes into communities defined as the maximal connected components. We denote communities by $\mathcal{C}_c$ with $c=1,\ldots,C$, where $C$ is the total number of communities. Each community has exactly one cycle. The length $\ell_c$ of the cycle in community $\mathcal{C}_c$ satisfies the obvious bounds: $1\leq \ell_c\leq |\mathcal{C}_c|$. The number of elements in the set $S$ is $M = |S| = \sum_{1\leq c\leq C} |\mathcal{C}_c|$.

Let us interpret elements of set $S$, equipped with map $f: S\to S$, as individuals. An individual $y\in S$ is an expert if there is at least one individual $x\in S$ such that $f(x)=y$; followers are individuals who are not experts. Denote by $E_c$ and $F_c$ the numbers of experts and followers in the community $\mathcal{C}_c$. The total numbers of experts and followers are $E =\sum_{1\leq c\leq C} E_c$ and $F =\sum_{1\leq c\leq C} F_c$. Every community decomposes into experts and followers, $ E_c+ F_c=|\mathcal{C}_c|$, and similarly $E+F = |S| = M$. We define egocentrics as loops in the graph: $f(x)=x$ for an egocentric $x\in S$. Introverts are isolated egocentrics, i.e., egocentrics without followers. 

\begin{figure}
\centering
\includegraphics[width=7.89cm]{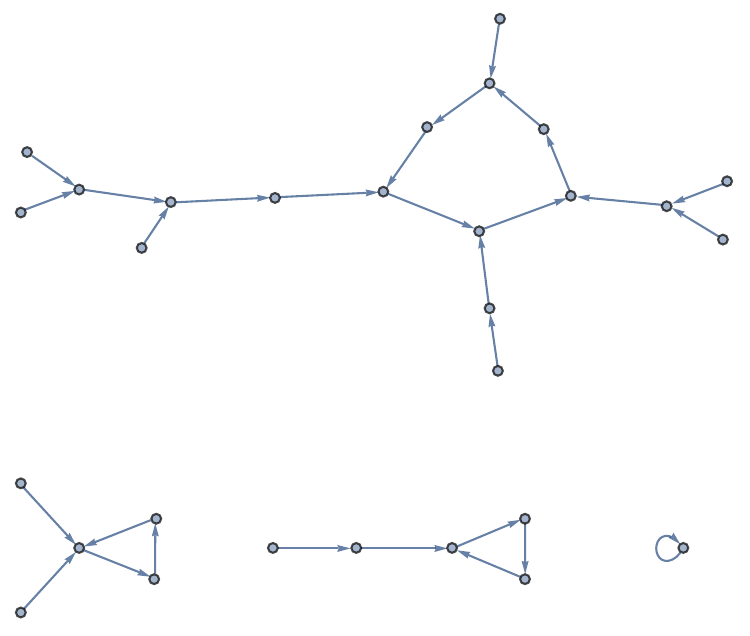}
\caption{An example of a random map with four communities. Top: The biggest community with $|\mathcal{C}_4|=18$ individuals and $F_4=7$ followers has a cycle of length $\ell_4=6$. Bottom: Communities $\mathcal{C}_3, \mathcal{C}_2, \mathcal{C}_1$ (from left to right). The first two communities have $|\mathcal{C}_3|=|\mathcal{C}_2|=5$ individuals, $F_3=2$ and $F_2=1$ followers, and cycles of length $\ell_3=\ell_2=3$. The community $\mathcal{C}_1$ consists of an introvert (isolated egocentric): A single individual $|\mathcal{C}_1|=1$ with no followers $F_1=0$.  }
\label{Fig:random-map-ill}
\end{figure}

In an example in Fig.~\ref{Fig:random-map-ill},  there are four communities with $|\mathcal{C}_1|=1, |\mathcal{C}_2|=|\mathcal{C}_3|=5, |\mathcal{C}_4|=18$  individuals and $F_1=0, F_2=1, F_3=2, F_4=7$ followers; overall there are 10 followers and 19 experts. The lengths of the cycles are $\ell_1=1, \ell_2=\ell_3=3, \ell_4=6$. There is one egocentric who happens to be an introvert. 

For a random map $f: S\to S$, the image  $f(x)$ of each $x\in S$ is  chosen randomly. For random maps without any constraint (classical random maps), all $M^M$ maps are considered equally likely. See \cite{Kruskal,Katz,Stepanov69,Stepanov71} for earlier work on classical random maps and \cite{Harris,Kolchin} for reviews of earlier research. Random maps have applications in modeling epidemic processes, analyses of cryptographic systems, and random number generation, etc. (see \cite{Gertsbakh,Pittel83,Flajolet90,Quisquater90,Jaworski08,Aldous,Jaworski14} and references therein). Random maps appear in the context of spin glasses and disordered materials \cite{Witten-RM}. Special types of random maps arise in the context of synchronous dynamics of random Boolean networks such as Kauffman's NK model \cite{BN-Kauffman69,BN-Kauffman89,Kauffman,BN-Hilhorst87,BN-Derrida87,BN-Flyvbjerg,BN-Parisi,BN-Sam,BN-Barbara-1,BN-Barbara-2,BN-Mori}.  Random maps can be used to mimic belief dynamics \cite{Mirta}. Other applications are anticipated in social and communication networks \cite{Wasserman,Castellano09,Tomasello}.

Even for classical random maps (Sec.~\ref{sec:maps}), some results related to followers, experts, egocentrics, and introverts are scattered in the literature \cite{Harris,Kolchin,Flajolet90}, e.g., terminal nodes appearing in Ref.~\cite{Flajolet90} are experts in our language. The sociological terminology can play an assertive role, suggesting novel classes of random maps mimicking social networks. In this paper, we propose and analyze random maps with sociological favor. Specifically, we introduce random maps compatible with fitness (Sec.~\ref{sec:FM}) and ranking (Sec.~\ref{sec:MRM}). 

\begin{table}[h!]
\centering
\renewcommand{\arraystretch}{1.5}{
\begin{tabular}{|c|c|c|c|c| c|}
\hline
Maps           & Fitness      & Uniform            & Linear                   & $\text{(Linear)}^{-1}$  \\
\hline
 Experts       & $e(a)$     & $\frac{1}{2}$    &$\frac{9-e^2}{4}$     &       $ 1-(\ln M)^{-1}$  \\
\hline
 Egocentrics & 1             & $\ln M$            &$2$     &       $ M/\ln M$  \\
\hline
 Introverts     & $A(a)$     & 1      &$(e^2-1)/M$ & $ e^{-1} M/\ln M$  \\
\hline
\end{tabular}
}
\caption{We give the fraction of experts and the average numbers of egocentrics and introverts for the class of fitness models  \eqref{fitness-a} with $a>-1$,  and for random map models [\eqref{map:more}, \eqref{map:linear-rank}, \eqref{map-linear} from left to right] compatible with ranking. For the general class of fitness models \eqref{map:fitness}, there is, on average, one egocentric individual independently of $M$.  } 
\label{Tab:loops}
\end{table}

In Table~\ref{Tab:loops}, we collect a few results for three random map models compatible with ranking and for a one-parametric class of fitness models in which we choose the image $j=f(i)$ for each $i\in S$ with probability proportional to $j^a$, the fitness of the image; see \eqref{fitness-a} for the precise definition. If $a>-1$, the fraction of experts $e(a)$ and the average number of introverts are (Sec.~\ref{sec:FM})
\begin{subequations}
\begin{align}
\label{e-a}
&e(a) = 1 - \int_0^1 dx\,e^{-(a+1)x^a}\\
\label{A-a}
&\langle I\rangle =  A(a) = (a+1)\int_0^1 dx\,x^a e^{-(a+1)x^a}
\end{align}
\end{subequations}

Random maps compatible with ranking (Sec.~\ref{sec:MRM}) are defined as follows. We label elements of the set $S$ from $1$ to $M$, interpret labels as ranks, and postulate that a map $f: S\to S$ is compatible with ranking if $f(i)\geq i$ for all $i=1,\ldots,M$. All $M!$ of maps compatible with ranking are considered equiprobable for recursive random map (RRM). For the RRM, the probability distribution of the random variable $F$ is 
\begin{equation}
\label{PFM:Euler}
\Pi(F,M)=\frac{1}{M!}\,\eulerian{M}{F}
\end{equation}
where $\eulerian{M}{F}$ are Eulerian numbers \cite{Euler36,Euler55,Knuth}. We also derive an exact formula for the probability distribution of the number of communities
\begin{equation}
\label{PCM:Stirling}
P(C,M)=\frac{1}{M!}\,{M\brack C}
\end{equation}
where ${M\brack C}$ are the Stirling numbers of the first kind \cite{Knuth}. The degree $k$ of the individual with highest rank, the founder, has identical probability distribution
\begin{equation}
\label{PkM:Stirling}
\Phi(k,M)=\frac{1}{M!}\,{M\brack k}
\end{equation}
By definition \cite{Knuth}, the Eulerian numbers $\eulerian{M}{F}$ and Stirling numbers ${M\brack C}$  count special types of permutations. These numbers are intimately related to the symmetric group $S_M$ of permutations of the numbers $1$ to $M$, and more generally to group theory (see, e.g., \cite{Eulerian80,Coxeter}).

Some static random maps admit a dynamical reformulation that may be more amenable to analysis. Some features make no sense in the static framework but naturally arise in the dynamical realm. The prime static random map admitting a dynamical reformulation is the RRM model. Indeed, the RRM is isomorphic to a growing random map where individuals arrive one by one, and each new individual maps to an individual randomly chosen amongst those already present (including itself). Stopping the process after $M$ steps gives the RRM of size $M$. In Appendix~\ref{ap:RRM}, we explain this isomorphism and derive \eqref{PFM:Euler} and \eqref{PCM:Stirling}. Other evolving mechanisms are possible, and in Sec.~\ref{sec:rewire}, we briefly discuss random maps of the same set evolving through rewiring. In Sec.~\ref{sec:prophet}, we define prophets and founders and briefly analyze their properties in random map models. In particular, we derive the degree distribution \eqref{PkM:Stirling} of the founder. We conclude with a discussion (Sec.~\ref{sec:disc}).

\section{Classical Random Maps}
\label{sec:maps}

In this section, we consider classical random maps: all $M^M$ maps of a set of size $|S|=M$ to itself are equiprobable. Many results in this and the following sections are valid for arbitrary $M$, but we are primarily interested in the behavior in the $M\to\infty$ limit. 

We begin by outlining basic statistical properties of classical random maps. Denote by $P(C,M)$ the probability that there are $C$ communities. The average number of communities, $\langle C(M)\rangle = \sum_{C\geq 1}C P(C,M)$, grows logarithmically \cite{Kruskal,Harris} with $M$. More precisely \cite{Flajolet}
\begin{equation}
\label{c-av}
\langle C(M)\rangle = \tfrac{1}{2}\ln M + \tfrac{1}{2}(\gamma+\ln 2) + O(M^{-1/2})
\end{equation}
where $\gamma \approx 0.5772156649$ is the Euler constant.

The full probability distribution $P(C,M)$ is also known  \cite{Harris}. The probability to have exactly one community is \cite{Katz}
\begin{equation}
\label{P1-M}
P(1,M) = \frac{(M-1)!}{M^M}\sum_{m=0}^{M-1} \frac{M^m}{m!}
\end{equation}
For large $M$, the probability $P(1,M)$ simplifies to
\begin{equation}
\label{P1}
P(1,M)\simeq \sqrt{\frac{\pi}{2M}}
\end{equation}

The probability $\Pi_1(\ell,M)$ to end up with a single community that has a cycle of length $\ell$ is also known \cite{Katz}:
\begin{equation}
\Pi_1(\ell,M) = \frac{M!}{(M-\ell)! M^{\ell + 1}}
\end{equation}
This probability distribution simplifies to
\begin{equation}
\Pi_1(\ell,M) = M^{-1}\,e^{-x^2/2}
\end{equation}
when $M\to\infty$ and $\ell\to\infty$  with $x = \ell/\sqrt{M}$ kept finite.

A randomly chosen individual is a follower with probability $(1-1/M)^M$. Thus 
\begin{equation}
\label{Fav:random}
\langle F\rangle = M(1-1/M)^M
\end{equation}
For classical random maps, followers and experts were called terminal nodes and image points in \cite{Flajolet90} where the asymptotic behaviors $\langle F\rangle = Me^{-1}$ and $\langle E\rangle = M(1-e^{-1})$ were derived.

The probability distribution $Q(\mathcal{L},M)$ to have $\mathcal{L}$ loops is the binomial distribution 
\begin{equation}
\label{Q:binom}
Q(\mathcal{L},M) = \binom{M}{\mathcal{L}}\left(1-\frac{1}{M}\right)^{M-\mathcal{L}} \left(\frac{1}{M}\right)^\mathcal{L}
\end{equation}
The average number of loops is independent on $M$
\begin{equation}
\label{loops-av}
\langle \mathcal{L}\rangle = \sum_{\mathcal{L}=0}^M \mathcal{L} Q(\mathcal{L},M) = 1
\end{equation}
while higher moments depend on $M$:
\begin{equation}
\label{loops-higher}
\langle \mathcal{L}^2\rangle = 2 - M^{-1}, \quad \langle \mathcal{L}^3\rangle = 5 - 6M^{-1}+ 2 M^{-2}
\end{equation}
etc. Specializing \eqref{Q:binom} to $\mathcal{L}=0$ and $\mathcal{L}=M$ we otain the probabilities to observe extremal numbers of loops:
\begin{equation}
\label{Q:min-max}
Q(0,M) =  \big(1-M^{-1}\big)^M,  \quad Q(M,M) =  M^{-M}
\end{equation}
The probability to have maximal number of communities, $C=M$, is of course  identical to the probability to have the maximal number of loops:  
\begin{equation}
\label{QMM:classical}
P(M,M) = Q(M,M) =  M^{-M}
\end{equation}

In the $M\to \infty$ limit, the binomial distribution \eqref{Q:binom} simplifies to the Poisson distribution
\begin{equation}
\label{Q:Poisson}
Q(\mathcal{L},\infty)  = \frac{e^{-1}}{\mathcal{L}!}
\end{equation}

A loop $f(i)=i$ implies that individual $i$ is egocentric. Thus $Q(\mathcal{L},M)$ is the probability that there are $\mathcal{L}$ egocentrics. An introvert is an isolated egocentric, i.e., someone who is not an expert for anyone else. The probability $R(I,M)$ to have $I$ introverts can be found recurrently. The probability $R(1,M)$ satisfies 
\begin{equation}
\label{R1}
R(1,M) = M^{-1} \binom{M}{1} \left(1-\frac{1}{M}\right)^{M-1}R(0,M-1)
\end{equation}

Here $M^{-1}$ is the probability to form a loop, the binomial factor accounts for the label of an introvert, $(1-M^{-1})^{M-1}$ assures that remaining individuals map into themselves and with probability $R(0,M-1)$ there are no introverts among them. Similarly to \eqref{R1}, the probability $R(2,M)$ to have two introverts satisfies 
\begin{equation}
\label{R2}
R(2,M) = M^{-2}\binom{M}{2}\left(1-\frac{2}{M}\right)^{M-2}R(0,M-2)
\end{equation}

In the $M\to\infty$ limit, Eqs.~\eqref{R1}--\eqref{R2} become
\begin{equation}
\label{R12:inf}
R(1,\infty) = e^{-1}R(0,\infty),\quad R(2,\infty) = \frac{e^{-2}}{2!}\,R(0,\infty)
\end{equation}

Generally for arbitrary $I\geq 1$ we obtain 
\begin{equation}
\label{RI:inf}
R(I,\infty) = \frac{e^{-I}}{I!}\,R(0,\infty)
\end{equation}
The normalization requirement $\sum_{I\geq 0}R(I,\infty) = 1$ fixes the probability to have no introverts 
\begin{equation}
\label{R0}
R(0,\infty)  = \exp\!\big[-e^{-1}\big] = 0.6922006275553\ldots
\end{equation}
Thus $R(I,\infty)$ is the Poisson distribution
\begin{equation}
\label{RI:dist}
R(I,\infty) = \frac{1}{I!}\,\exp\!\big[-I -e^{-1}\big]
\end{equation}
Using \eqref{RI:dist} we find the average number of introverts
\begin{equation}
\label{introverts-av}
\langle I\rangle = e^{-1}=0.3678794411714\ldots
\end{equation}
Equations \eqref{loops-av} and \eqref{introverts-av} give average numbers, not fractions---in classical random maps, the numbers of egocentrics and introverts remain finite when $M\to\infty$. 

The degree of the vertex is the number of its neighbors. Followers have degree 1, while experts have degrees larger than one. It makes sense to stratify experts according to their degree $k$. Denote by $N_k(M)$ the number of individuals of degree $k$. We have
\begin{equation}
\label{Nk:EF}
F = N_1, \quad E = \sum_{k\geq 2} N_k(M)
\end{equation}

The notion of degree is easy to appreciate in undirected graphs. In a directed graph, each vertex has two natural degrees: an in-degree $k_\text{in}$ and an out-degree $k_\text{out}$. The full degree of the vertex is the sum
\begin{equation}
k_\text{in}  + k_\text{out} = k
\end{equation}
For directed graphs generated by maps all out-degrees are equal to unity, $k_\text{out}(x)=1$ for each $x\in S$. Therefore it suffices to consider the full degree. If $k(x)=1$, vertex $x$ is a follower; if $k(x)\geq 2$, vertex $x$ is an expert. Each edge contributes a unity to the total in and out degrees, 
\begin{equation*}
\sum_{x\in S}k_\text{in}(x) = \sum_{x\in S}k_\text{out}(x) = M, 
\end{equation*}
and therefore
\begin{equation}
\label{sum:deg}
\sum_{x\in S}k(x) = 2M 
\end{equation}
for an arbitrary map. The sum rules
\begin{subequations}
\begin{align}
\label{norm}
& \sum_{k\geq 1} N_k(M) = M\\
\label{map:def}
& \sum_{k\geq 1} kN_k(M) = 2M
\end{align}
\end{subequations}
are valid for an arbitrary map: \eqref{norm} follows from \eqref{Nk:EF} and $E+F=M$, while \eqref{map:def} follows from \eqref{sum:deg}. 

For classical random maps, and for most models studied below, $N_k(M)$ are random self-averaging quantities. Therefore, the fractions 
\begin{equation}
n_k(M) = M^{-1} \langle N_k(M) \rangle
\end{equation}
saturate in the $M\to\infty$ limit and hence provide the major information about the asymptotic behavior of $N_k(M)$. From the definition of classical random maps, we find 
\begin{equation}
\label{nk:binom}
n_k(M) = \binom{M}{k-1}\left(1-\frac{1}{M}\right)^{M-k+1} \left(\frac{1}{M}\right)^{k-1}
\end{equation}
for arbitrary $M$. The binomial distribution \eqref{nk:binom} becomes the Poisson distribution 
\begin{equation}
\label{nk:Poisson}
n_k(\infty) = \frac{e^{-1}}{(k-1)!}
\end{equation}
when $M\to\infty$. This result appears in \cite{Flajolet90}. 

The knowledge of the degree distribution allows one to probe the asymptotic behavior of the maximal degree $D(M)$. Using the criterion 
\begin{equation}
\label{D:criterion}
M\sum_{k\geq D(M)}n_k\sim 1
\end{equation}
and \eqref{nk:Poisson} we find  
\begin{equation}
\label{DM:nu}
D(M) \simeq  \frac{\ln M}{\ln(\ln M)}  
\end{equation}

\begin{figure}
\centering
\includegraphics[width=7.77cm]{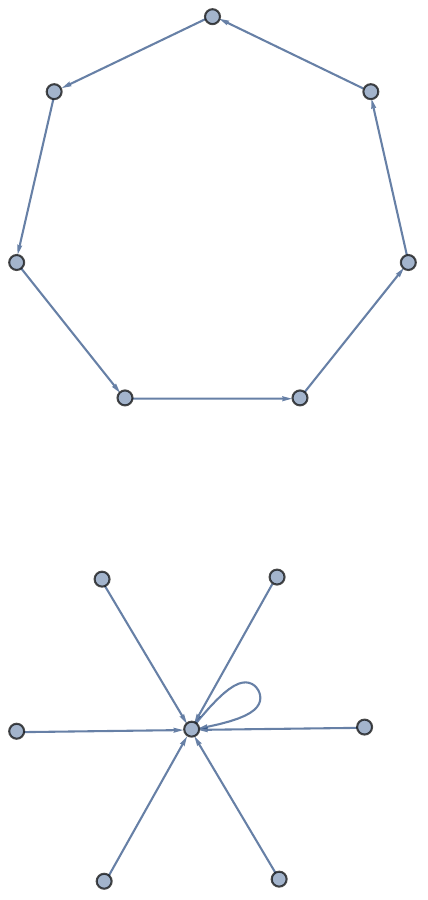}
\caption{The community on the top is a cycle. All individuals in a cycle are experts. The number of followers is minimal, $F=0$, for any collection of cycles. The graph on the bottom is a directed star with one expert.  The number of followers is maximal, $F=M-1$, if the map is a directed star. }
\label{Fig:graphs}
\end{figure}

Finally we mention the bounds
\begin{equation}
\label{F:bounds}
0\leq F\leq M-1
\end{equation}
for total number of followers. The lower bound is achieved when all communities are cycles (Fig.~\ref{Fig:graphs}). Each such map is a permutation, so the number of such maps is $M!$, and hence the probability to have no followers is
\begin{equation}
\label{F-min:RM}
\Pi(0,M) = \frac{M!}{M^M}
\end{equation}
The upper bound in \eqref{F:bounds} is achieved when a map is a directed star, see Fig.~\ref{Fig:graphs}. The probability to have the maximal number of followers is
\begin{equation}
\label{F-max:RM}
\Pi(M-1,M) = \frac{1}{M^{M-1}} 
\end{equation}

\section{Fitness Models}
\label{sec:FM}

In this section, we introduce and investigate random map models compatible with fitness. We assign a positive number, fitness, to each individual. We label individuals from 1 to $M$, denote by $\phi(i)$ the fitness of $i\in S$, and postulate that the probability of a map $i\to j$ is proportional to the fitness $\phi(j)$ of the image. Thus
\begin{equation}
\label{map:fitness}
\text{Prob}[f(i)=j] = \frac{\phi(j)}{Z}\,, \quad Z = \sum_{k=1}^M \phi(k)
\end{equation}
Classical random maps correspond to the uniform fitness:  $\phi(i)=1$ for all $i\in S$. 

More generally, one can postulate that the probability of a map $i\to j$ is proportional to $\phi_i(j)$, the fitness of $j$ in eye of $i$. Thus 
\begin{equation}
\label{gen:fitness}
\text{Prob}[f(i)=j] = \frac{\phi_i(j)}{Z_i}, \quad Z_i = \sum_{k=1}^M \phi_i(k)
\end{equation}
Each model \eqref{gen:fitness} is specified by $M^2$ fitness parameters $\phi_i(j)$. (Multiplying all $\phi_i(j)$ by the same factor does not change the rule \eqref{gen:fitness}, so $M^2-1$ parameters suffice.)

In this section,  we limit ourselves to the class of models \eqref{map:fitness} if not stated otherwise. Vertex $i$ is a follower with probability $[1-Z^{-1}\phi(i)]^M$, so 
\begin{equation}
\label{Lav:fitness}
\langle F\rangle = \sum_{i=1}^M \left[1-\frac{\phi(i)}{Z}\right]^M 
\end{equation}
is the average number of followers. To deduce concrete results from the general solution \eqref{Lav:fitness}, we consider special fitness models. For  the linear fitness model, $\phi(i)=i$, the large $M$ behavior of the sum in \eqref{Lav:fitness} reads
\begin{equation*}
\begin{split}
 &\left[1-\frac{\phi(i)}{Z}\right]^M  =  \left[1-\frac{2i}{M(M+1)}\right]^M \simeq e^{-2i/M}\\
 & \sum_{i=1}^M e^{-2i/M} = \frac{1-e^{-2}}{e^{2/M}-1} \simeq  \frac{1-e^{-2}}{2}\,M
 \end{split}
\end{equation*}
Thus
\begin{equation}
\label{L-av-fitness}
\lim_{M\to\infty}\frac{\langle E\rangle}{M} = \frac{1+e^{-2}}{2}\,, \quad   \lim_{M\to\infty}\frac{\langle F\rangle}{M} = \frac{1-e^{-2}}{2}
\end{equation}

More generally, consider a one-parameter family of fitness models with algebraic fitness $\phi(i)=i^a$, that is
\begin{equation}
\label{fitness-a}
\text{Prob}[f(i)=j] = \frac{j^a}{Z}\,, \quad Z = \sum_{k=1}^M k^a
\end{equation}
Interesting behaviors occur when $a\geq -1$. In the following, we tacitly assume that $a>-1$ if not stated otherwise; occasionally, we comment on subtle behaviors in the marginal $a=-1$ case. 

For the one-parameter family  \eqref{fitness-a} of fitness models, the average number of followers \eqref{Lav:fitness} becomes
\begin{equation*}
\langle F\rangle = \sum_{i=1}^M \left[1-\frac{i^a}{Z}\right]^M, \quad Z=\sum_{i=1}^M i^a
\end{equation*}
When $a>-1$, the `partition function' $Z$ diverges as $M^{a+1}/(a+1)$ when $M\gg 1$. This simplifies the asymptotic analysis: One replaces summation by integration and arrives at the integral representation
\begin{equation}
\label{Fav}
f(a) = \int_0^1 dx\, e^{-(a+1)x^a} 
\end{equation}
of the asymptotic fraction $f=\lim_{M\to\infty}M^{-1}\langle F\rangle$ of the followers. The integral representation \eqref{Fav} is applicable for all $a>-1$. When $a>0$, one can reduce \eqref{Fav} to
\begin{equation*}
f(a) = \frac{\Gamma(a^{-1})-\Gamma(a^{-1}, a+1)}{a(a+1)^{1/a}}
\end{equation*}
containing an incomplete gamma function
\begin{equation*}
\Gamma(\beta,B)=\int_B^\infty dy\,y^{\beta-1} e^{-y}
\end{equation*}
The fraction of experts is complementary as $E+F=M$; this fraction is given by \eqref{e-a}. 

\begin{figure}
\centering
\includegraphics[width=7.77cm]{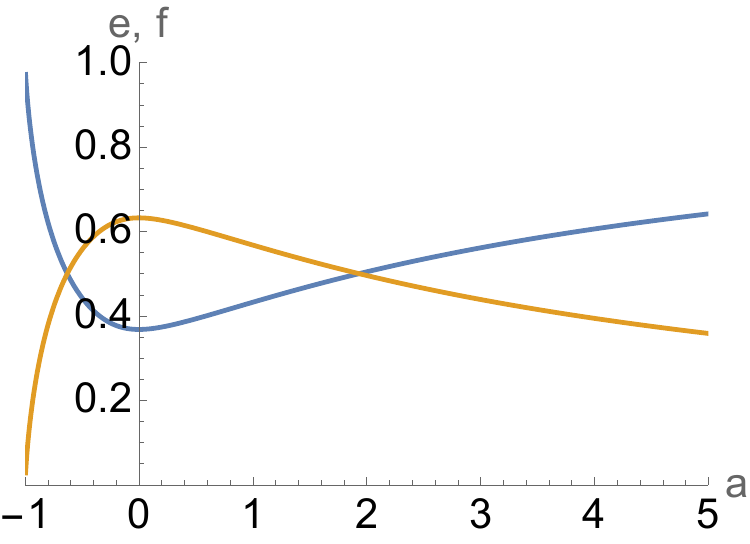}
\caption{The fractions of experts, Eq.~\eqref{e-a}, and followers, Eq.~\eqref{Fav}, for the family of fitness models \eqref{fitness-a}. The fraction of experts reaches maximum $1-e^{-1}$ for classical random maps ($a=0$),  and it vanishes when $a\to -1$ or $a\to\infty$. The fraction of followers exhibits a complementary behavior. }
\label{Fig:followers_experts}
\end{figure}

The dependence of the fractions $e(a)$ and $f(a)$ on the parameter $a$ is presented in Fig.~\ref{Fig:followers_experts}. For classical random maps, $a=0$, the fraction of experts is maximal while the fraction of followers is minimal among all models \eqref{fitness-a}. If $a\to\infty$ (we first take the $M\to\infty$ limit), the fraction of experts approaches zero as 
\begin{equation*}
\label{E-av:asymp}
e(a) = \frac{\ln a + \gamma}{a} -\frac{(\ln a)^2+2\gamma \ln a+\frac{\pi^2}{6}+\gamma^2-2}{2 a^2}+ \ldots
\end{equation*}
At two special values of $a$ which are the roots of
\begin{equation}
\frac{1}{2} = \int_0^1 dx\, e^{-(a+1)x^a} 
\end{equation}
the fractions of experts and followers are equal; this occurs when $a\approx 1.9329$ and $a\approx -0.6289$.  

The behavior is particularly subtle when $a=-1$. We need to compute
\begin{equation}
\label{Eav:1}
\langle E\rangle = \sum_{k=1}^M \left(1-\left[1-\frac{1}{k H_M}\right]^M\right)
\end{equation}
where $H_M= \sum_{1\leq m\leq M} m^{-1}$ is the harmonic number. Harmonic numbers $H_M$ grow logarithmically with $M$. More precisely 
\begin{equation}
\label{Harmonic}
H_M = \ln M + \gamma + \frac{1}{2M} -  \frac{1}{12M^2} + \frac{1}{120M^4} + \ldots 
\end{equation}
when $M\gg 1$. 

Replacing summation over $1\leq k\leq M$ in \eqref{Eav:1} by integration over $\kappa=k H_M/M$ we reduce \eqref{Eav:1} to
\begin{equation}
\label{Eav:int}
\frac{\langle E\rangle}{M} \simeq \frac{1}{H_M}\int_{M^{-1}H_M}^{H_M}d\kappa \,\big(1-e^{-1/\kappa}\big)
\end{equation}
The dominant contribution to the integral is gathered when $\kappa$ is large and equals to $\ln H_M$ in the leading order. Using the asymptotic  \eqref{Harmonic} we arrive at
\begin{equation}
\label{Eav:asymp}
\frac{\langle E\rangle}{M} \simeq \nu, \qquad \nu \simeq \frac{\ln(\ln M)}{\ln M}
\end{equation}
Therefore, when $a=-1$, the fraction of experts vanishes anomalously slowly in the $M\to\infty$ limit. The above asymptotic formula gives the leading behavior of $\nu$ when $\ln(\ln M)\gg 1$, i.e., for astronomically large $M$. A more careful estimate of the integral in \eqref{Eav:int}  and together with the asymptotic \eqref{Harmonic} lead to a more accurate asymptotic behavior of $\nu(M)$:
\begin{equation*}
\nu \simeq \frac{\ln(\ln M)+1-\gamma}{\ln M+\gamma}
\end{equation*}

We now turn to egocentrics, i.e., we investigate loops. For the class of models \eqref{map:fitness},  vertex $i$ is the loop with probability $Z^{-1}\phi(i)$.  Thus on average, there is one egocentric individual, no dependence on fitness and the population size $M$:
\begin{equation}
\langle \mathcal{L}\rangle = \sum_{i=1}^M \frac{\phi(i)}{Z} = 1
\end{equation}

To probe the probability distribution $Q(\mathcal{L},M)$, we first notice that there are no loops with probability 
\begin{equation}
\label{Q0M}
Q(0,M) = \prod_{i=1}^M \left[1-\frac{\phi(i)}{Z}\right]
\end{equation}
The probability \eqref{Q0M} tends to $Q(0,\infty)  = e^{-1}$ for all fitness models \eqref{map:fitness} in which $Z$ diverges when $M\to\infty$. For instance, for the class of models \eqref{fitness-a} with $a\geq -1$, the asymptotic behavior of $Q(0,M)$ is
\begin{equation}
\label{Q0M:asymp}
\frac{Q(0,M)}{Q(0,\infty)} - 1 \simeq - 
\begin{cases}
\frac{(a+1)^2}{2(2a+1)}\,M^{-1} &  a>-1\\
\frac{\pi^2}{12}\,(\ln M)^{-2}       &  a=-1
\end{cases}
\end{equation}

The probability to have one loop is
\begin{equation}
Q(1,M) = \prod_{i=1}^M \left[1-\frac{\phi(i)}{Z}\right] \sum_{j=1}^M \frac{\phi(j)}{Z-\phi(j)}
\end{equation}
If $Z$ diverges in the $M\to\infty$ limit, $Q(1,M)$ approaches to the same universal value $Q(1,\infty)  = e^{-1}$ as the probability to have no loops.

Generally fixing $\mathcal{L}$ and taking the $M\to\infty$ limit, one finds that $Q(\mathcal{L},\infty)$ is the Poisson distribution \eqref{Q:Poisson} as in the case of random maps. It is remarkable that $Q(\mathcal{L},\infty)$ is the same for all models \eqref{map:fitness} modulo just one property: $Z$ must diverge in the $M\to\infty$ limit. 

The dependence of $Q(\mathcal{L},M)$ on fitnesses becomes appreciable when $\mathcal{L}$ is large. For instance, the number of loops is maximal with probability
\begin{equation}
\label{QMM:fitness}
Q(M,M) = \prod_{i=1}^M \frac{\phi(i)}{Z}
\end{equation}

For models \eqref{fitness-a} with small integer values of the exponent $a$ the number of loops is maximal with probability
\begin{equation}
\label{QMM:0123}
Q(M,M) = 
\begin{cases}
M^{-M}                                                              & a=0\\
\frac{2^M \cdot M!}{[M(M+1)]^M}                      & a=1\\
\frac{6^M \cdot (M!)^2}{[M(M+1)(2M+1)]^M}    & a=2\\
\frac{4^M \cdot (M!)^3}{[M(M+1)]^{2M}}            & a=3
\end{cases}
\end{equation}

For the class of models \eqref{map:fitness}, the probability to have the maximal number of followers is
\begin{equation}
\label{F-max:fitness}
\Pi(M-1,M) = \sum_{i=1}^M \left[\frac{\phi(i)}{Z}\right]^M
\end{equation}
For classical random maps, $\phi(i)=1$, we recover \eqref{F-max:RM}. For the model with linear fitness, $\phi(i)=i$, we get
\begin{equation}
\Pi(M-1,M) = \frac{2^M}{[M(M+1)]^M} \,\sum_{i=1}^M i^M 
\end{equation}
For models \eqref{fitness-a}, the asymptotic behavior is
\begin{equation}
\label{F-max:a-asymp}
\Pi(M-1,M) \simeq e^{-\frac{a+1}{2}}\left(\frac{a+1}{M}\right)^M
\end{equation}
 when $a>-1$, and $\simeq (H_M)^{-M}$ when $a=-1$.

Let us now look at the distribution of introverts. Extending the argument that led to \eqref{R1} we obtain
\begin{equation}
\label{R1:fitness}
\frac{R(1,M)}{R(0,M-1)} = \sum_{i=1}^M \frac{\phi(i)}{Z} \left[1-\frac{\phi(i)}{Z}\right]^{M-1}
\end{equation}
Similarly
\begin{equation}
\label{R2:fitness}
\frac{R(2,M)}{R(0,M-2)} = \sum_{1\leq i<j\leq M} \varphi(i) \, \varphi(j) 
\end{equation}
where we have used the shorthand notation
\begin{equation}
\varphi(i) = \frac{\phi(i)}{Z} \left[1-\frac{\phi(i)}{Z}\right]^{M-1}
\end{equation}
Generally
\begin{equation}
\label{RI:fitness}
\frac{R(I,M)}{R(0,M-I)} = \sum_{1\leq i_1<\ldots<i_I\leq M} \prod_{p=1}^I\varphi(i_p)
\end{equation}

\begin{figure}
\centering
\includegraphics[width=7.77cm]{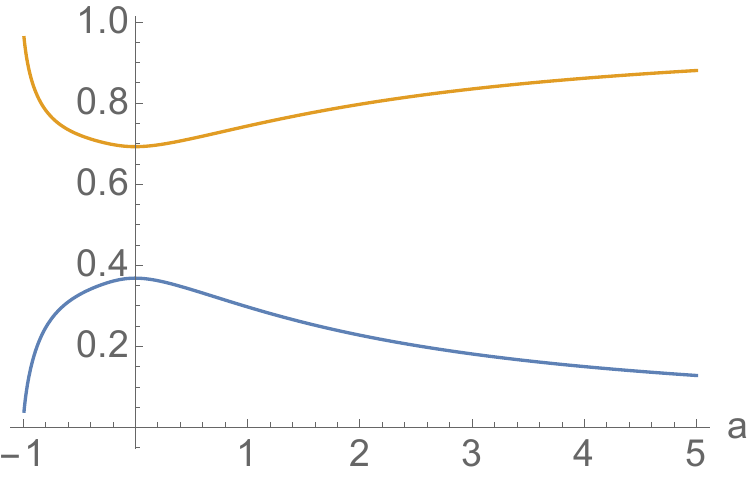}
\caption{Bottom: The average number of introverts, $A(a)$, see \eqref{R10} and \eqref{introverts-av-a}. Top: The probability $e^{-A(a)}$ that there are no introverts. Shown is the dependence on $a$ parametrizing fitness models \eqref{fitness-a}. The average number $A(a)$ of introverts coincides with the fraction of vertices of degree 2, see \eqref{n2AI}. }
\label{Fig:introverts}
\end{figure}

To deduce more concrete results from \eqref{RI:fitness} we consider again the one-parameter family of fitness models \eqref{fitness-a}. Taking the $M\to\infty$ limit and replacing summation in Eq.~\eqref{R1:fitness} by integration we obtain
\begin{equation}
\label{R10}
\frac{R(1,\infty)}{R(0,\infty)} = A(a) \equiv (a+1) \int_0^1 dx\, x^a e^{-(a+1)x^a}
\end{equation}
when $a>-1$. Similarly Eq.~\eqref{R2:fitness} simplifies to
\begin{equation}
\label{R20}
\frac{R(2,\infty)}{R(0,\infty)} = \frac{A^2}{2!}
\end{equation}
Continuing we arrive at the Poisson distribution
\begin{equation}
R(I,\infty) = \frac{A^I}{I!}\,e^{-A}
\end{equation}
The average number of introverts 
\begin{equation}
\label{introverts-av-a}
\langle I\rangle = A(a)
\end{equation}
and the probability $R(0,\infty) = e^{-A(a)}$ that there are no introverts are plotted in Fig.~\ref{Fig:introverts}.  The average number of introverts attains the maximal value at $a=0$, i.e. for the classical random map; $A(a)$ vanishes when $a\to -1$ and $a\to\infty$, e.g. $A(a)\simeq a^{-1}$ when $a\gg 1$. A few concrete values of $\langle I\rangle = A(a)$:
\begin{equation*}
\begin{split}
&A(0)=e^{-1}=0.3678794411714\ldots\\
&A(1)=\frac{1}{2}-\frac{3}{2e^2}=0.296997075145\ldots\\
&A(2)=\sqrt{\frac{\pi}{48}}\,\,\text{Erf}\!\left[\sqrt{3}\right]-\frac{1}{2e^3} =0.22727824593\ldots
\end{split}
\end{equation*}

The above results \eqref{R10}--\eqref{introverts-av-a} are valid when $a>-1$. A more subtle behavior occurs when $a=-1$. Strictly speaking, all $R(I,\infty)$ vanish for all $I\geq 1$, but the large $M$ behavior is still essentially described by the Poisson distribution 
\begin{equation}
R(I,M) = \frac{\nu^I}{I!}\,e^{-\nu}\,, \qquad \nu = \frac{\ln(\ln M)}{\ln M}
\end{equation}

We now turn to the degree distribution. For the class of models \eqref{map:fitness}, vertex $i$ has degree $k$ with probability
\begin{equation*}
\binom{M}{k-1}\left[1-\frac{\phi(i)}{Z}\right]^{M-k+1}\left[\frac{\phi(i)}{Z}\right]^{k-1}
\end{equation*}
Therefore the fraction of vertices with degree $k$ is
\begin{equation}
\label{nk:sum-fitness}
n_k = \sum_{i=1}^M \frac{1}{M}\binom{M}{k-1}\left[1-\frac{\phi(i)}{Z}\right]^{M-k+1}\left[\frac{\phi(i)}{Z}\right]^{k-1}
\end{equation}
Comparing \eqref{R1:fitness} and \eqref{nk:sum-fitness}, we deduce a curious result that the fraction of vertices of degree 2 is related to the distribution of introverts:
\begin{equation}
\label{R1n2}
n_2 = \frac{R(1,M)}{R(0,M-1)} 
\end{equation}
 For the one-parameter family of fitness models \eqref{fitness-a}, we use \eqref{R10} and \eqref{introverts-av-a} to obtain
\begin{equation}
\label{n2AI}
n_2 = A(a) = \langle I\rangle 
\end{equation}
when $M\to\infty$. 

\begin{figure}
\centering
\includegraphics[width=7.77cm]{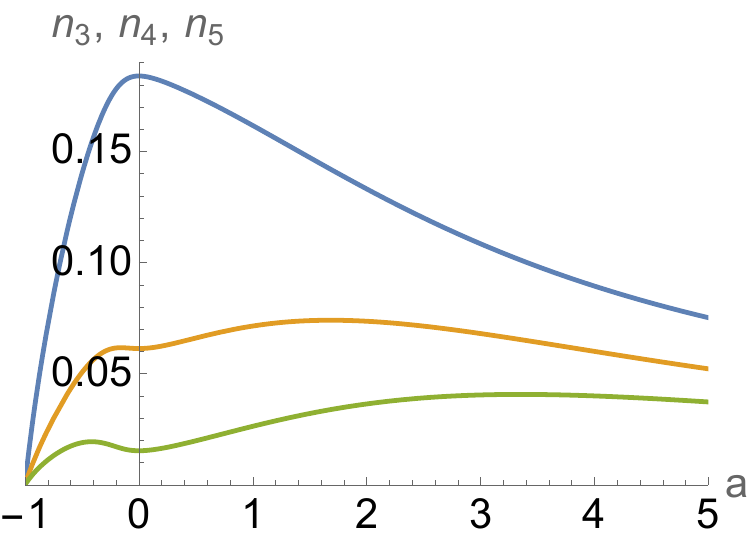}
\caption{The fitness model \eqref{map:fitness} with $\phi(i)=i^a$. The fractions $n_k$ with degrees $k=3,4,5$ (top to bottom) }
\label{Fig:n345}
\end{figure}

Specializing \eqref{nk:sum-fitness} to fitness models \eqref{fitness-a} and taking the $M\to\infty$ limit we obtain
\begin{equation}
\label{nk:fitness}
n_k = \frac{1}{(k-1)!}\int_0^1 dx\, \left[(a+1)x^a\right]^{k-1} e^{-(a+1)x^a}
\end{equation}
Setting $a=0$, we recover the Poisson degree distribution \eqref{nk:Poisson} describing classical random maps. The dependence of $n_1=M^{-1}\langle F\rangle$ on $a$ is shown in Fig.~\ref{Fig:followers_experts}. The dependence of $n_2=A(a)=\langle I\rangle$ on $a$ is shown in Fig.~\ref{Fig:introverts} (bottom curve). The fractions of vertices with degrees $k=3,4,5$ are shown in Fig.~\ref{Fig:n345}. 

The tail of the distribution \eqref{nk:fitness} changes from factorial to algebraic when $a$ passes through $a=0$ corresponding to random maps:
\begin{equation}
\label{nk:tail}
n_k \simeq 
\begin{cases}
a^{-1}e^{-a-1}\,\frac{(a+1)^{k-1}-1}{(k-1)(k-1)!} & a>0\\
- a^{-1}(1+a)^{-1/a}\,k^{-1+1/a}         &-1<a<0
\end{cases}
\end{equation}
when $k\gg 1$. The asymptotic behavior in the $a>0$ range is written in such a way that  it reduces to \eqref{nk:Poisson} in the $a\to +0$ limit. 

Equations \eqref{nk:fitness}--\eqref{nk:tail} apply when $a>-1$. In the marginal $a=-1$ case, all $n_k$ with $k\geq 2$ vanish when $M=\infty$. Keeping $M\gg 1$ large but finite, one gets
\begin{equation}
n_k \simeq
\begin{cases}
1 - \nu                                 & k=1\\
\nu                                      & k=2\\
\frac{(\ln M)^{-1}}{(k-2)(k-1)}     & k\geq 3
\end{cases}
\end{equation}
with $\nu$ defined in Eq.~\eqref{Eav:asymp}. 

Finally, we give an estimate of the maximal degree:
\begin{equation}
\label{DM:nu-a}
D(M) \simeq  
\begin{cases}
\nu^{-1}                   &  a\geq 0\\
\frac{a+1}{M^a}      & -1<a<0 \\
\frac{M}{\ln M}        & a= -1
\end{cases} 
\end{equation}
The remarkably universal growth law of $D(M)$ in the $a\geq 0$ range is established using the criterion \eqref{D:criterion} and the asymptotic \eqref{nk:tail}. A more precise estimate coming from \eqref{nk:tail} is implicitly determined by 
\begin{equation}
D(M)[\ln D(M) - 1 - \ln(a+1)]\simeq \ln M
\end{equation}
When $-1<a<0$, the maximal degree can be determined by first noting that the average degree of the vertex with label $i$ is $d_i=(a+1)i^a/M^a$. These average degrees are significantly separated (a gap between adjacent degrees greatly exceeds fluctuations of those degrees). Hence the maximal degree is $D(M)=d_1=(a+1)/M^a$ as stated in \eqref{DM:nu-a}. Similarly, when $a=-1$, the maximal degree is very close to the average degree $d_1=M/\ln M$ of the vertex with label 1.

\section{Random Maps Compatible with Ranking}
\label{sec:MRM}

Every finite set can be linearly ordered. Let us identify the ordering with labeling and interpret labels as ranks: $i<j$ asserts that individual $i$ has a lower rank than individual $j$. It is natural to consider maps compatible with ranking, namely maps $f:S\to S$ satisfying $f(i)\geq i$ for all elements $i\in S$. 

A subclass of the general fitness models \eqref{gen:fitness} satisfying $\phi_i(j) =\psi(j-i)$ and $\psi(n) = 0$ when $n<0$ is compatible with ranking. The simplest uniform choice is $\psi(n) = 1$ when $n\geq 0$, so for each $i$, an expert is chosen uniformly at random among all $j\geq i$. Thus 
\begin{equation}
\label{map:more}
\text{Prob}[f(i)=j] = 
\begin{cases}
0                          & 1\leq j<i\\
(M-i+1)^{-1}         &  i\leq j \leq M
\end{cases}
\end{equation}

If we want to mimic the desire to map to the higher rank individuals with higher probabilities, we may choose $\psi(n) = n+1$. In this model
\begin{equation}
\label{map:linear-rank}
\text{Prob}[f(i)=j] = 
\begin{cases}
0                                                 & 1\leq j<i\\
\frac{2(j-i+1)}{(M-i+1)(M-i+2)}     &  i\leq j \leq M
\end{cases}
\end{equation}

Choosing $\psi(n)$ decreasing with $n$ mimics peer affinity relevant in many social settings. For the inverse linear function $\psi(n) = 1/(n+1)$ 
\begin{equation}
\label{map-linear}
\text{Prob}[f(i)=j] = 
\begin{cases}
0                                                & 1\leq j<i\\
\frac{1}{(j-i+1) H_{M-i+1}}          &  i\leq j \leq M
\end{cases}
\end{equation}

A one-parametric family of ranking models interpolating between models \eqref{map:more}--\eqref{map-linear} is defined by the rule
\begin{equation}
\label{map-a}
\text{Prob}[f(i)=j] = 
\begin{cases}
0                                                & 1\leq j<i\\
\frac{(j-i+1)^a}{S_{M-i+1}(a)}     &  i\leq j \leq M
\end{cases}
\end{equation}
with $S_k(a)=\sum_{1\leq p\leq k} p^a$. The interesting range is again $a\geq -1$.  

For random maps compatible with ranking, cycles necessarily have unit length, i.e., each cycle is a loop. Hence, the probability of the minimal number of followers is equal to the probability to have the maximal number of loops:
\begin{equation}
\label{F-min:rank}
\Pi(0,M)=Q(M,M)
\end{equation}
The probability of the maximal number of followers is
\begin{equation}
\label{F-max:rank}
\Pi(M-1,M)=\prod_{i=1}^M \frac{(M-i+1)^a}{S_{M-i+1}(a)} 
\end{equation}
for the class of models \eqref{map-a}.

\subsection{Uniform choice: Model \eqref{map:more}}
\label{sec:rank1}

The average number of egocentric individuals, i.e., the average number of loops, is the harmonic number
\begin{equation}
\langle \mathcal{L}(M)\rangle = H_M = \sum_{m=1}^M \frac{1}{m}
\end{equation}
Indeed, the top rank individual maps to itself, the second highest rank individual maps to itself with probability 1/2,  the third highest rank individual maps to itself with probability 1/3, etc.  

There is always at least one loop as the top rank individual maps to itself, $M\to M$. The probability to have exactly one loop is
\begin{equation}
Q(1,M) = \prod_{i=2}^{M} \left[1-\frac{1}{i}\right] = \frac{1}{M}
\end{equation}
The probability to have exactly two loops is
\begin{equation}
Q(2,M) = \prod_{i=2}^{M} \left[1-\frac{1}{i}\right] \sum_{j=1}^{M-1} \frac{1}{j}= M^{-1}\, H_{M-1}
\end{equation}
The probability to have exactly three loops is
\begin{equation}
Q(3,M) =  \frac{1}{M} \sum_{1\leq j<k\leq M-1} \frac{1}{j k}
\end{equation}
Computing the sum one gets
\begin{equation}
Q(3,M) =  \frac{1}{2M} \left[(H_{M-1})^2-H_{M-1}^{(2)}\right]
\end{equation}
where $H_n^{(2)}=\sum_{1\leq k\leq n}k^{-2}$ is the generalized harmonic number. Thus 
\begin{equation}
Q(3,M) \simeq  \frac{1}{2M} \left[(\ln M + \gamma)^2-\frac{\pi^2}{6}\right]
\end{equation}
when $M\gg 1$. 

Generally, the  probability to have exactly $\mathcal{L}$ loops is given by 
\begin{equation}
\label{QLM}
Q(\mathcal{L},M) =  \frac{1}{M}  \sum_{1\leq j_1<\ldots<j_{\mathcal{L}-1}\leq M-1} \prod_{\nu=1}^{\mathcal{L}-1} \frac{1}{j_\nu}
\end{equation}
This formula is valid for all $2\leq \mathcal{L}\leq M$. When $\mathcal{L}$ is finite and $M\gg 1$,  Eq.~\eqref{QLM} simplifies to 
\begin{equation}
Q(\mathcal{L},M)  \simeq  M^{-1}\,\frac{(\ln M)^{\mathcal{L}-1}}{(\mathcal{L}-1)!}
\end{equation}
If $\mathcal{L}=M$, the sum in Eq.~\eqref{QLM} consists of a single term, so the maximal number of loops occurs with probability
\begin{subequations}
\begin{equation}
\label{QMM:RRM}
Q(M,M) =  \frac{1}{M!}
\end{equation}
With a little more effort, one extracts
\begin{align}
\label{QMM-1:RRM}
Q(M-1,M) &=  \frac{M(M-1)}{2 M!}\\
\label{QMM-2:RRM}
Q(M-2,M) &=  \frac{M(M-1)(3M-1)}{24 M!}
\end{align}
\end{subequations}
from Eq.~\eqref{QLM}. 

To determine the average number of followers we note that vertex $1$ is a follower with probability  $1-1/M$ as it can be an expert only to itself; vertex $2$ is a follower with probability  $[1-\frac{1}{M}][1-\frac{1}{M-1}]=1-2/M$; generally vertex $i$ is a follower with probability $1-i/M$. Thus  
\begin{equation}
\label{F:MRM}
\langle F\rangle = \sum_{i=1}^M \left[1-\frac{i}{M}\right] = \frac{M-1}{2}
\end{equation}

In Appendix~\ref{ap:RRM}, we show that the model \eqref{map:more} is isomorphic to recursive random maps and derive Eq.~\eqref{PFM:Euler} for the probability distribution of the number of followers
\begin{equation}
\label{PFM:def}
\Pi(F,M) = \text{Prob}[F(M)=F] 
\end{equation}
We also derive the announced formula \eqref{PCM:Stirling} expressing the probability distribution of the number of communities 
\begin{equation}
\label{PCM:def}
P(C,M) = \text{Prob}[C(M)=C]
\end{equation}
via the Stirling numbers of the first kind \cite{Knuth}. Using \eqref{PCM:Stirling}, we express the average number of communities via the harmonic number
\begin{equation}
\label{C-av-RRM}
\langle C\rangle = \sum_{C\geq 1}C P(C,M) = H_M
\end{equation}
Thus, the average number of communities grows logarithmically with $M$ as for classical random maps,  Eq.~\eqref{c-av}, but with twice larger amplitude [cf. \eqref{Harmonic}].

The distributions $P(C,M)$ and $\Pi(F,M)$ do not become stationary (i.e., independent on $M$) in the $M\to\infty$ limit. This is obvious from explicit formulae for $P(C,M)$ at small $C$ such as
\begin{equation}
\label{P12}
P(1,M)=M^{-1}, \qquad P(2,M)=M^{-1}H_{M-1}
\end{equation}
and from explicit formulae for $\Pi(F,M)$ at small $F$
\begin{equation}
\label{P01}
\Pi(0,M)=\frac{1}{M!}, \qquad \Pi(1,M)=\frac{2^M-M-1}{M!}
\end{equation}

The distribution of introverts $R(I,M)$ is easier to compute than the distributions $P(C,M)$ and $\Pi(F,M)$. The key observation is that each individual is introvert with probability $1/M$. This assertion is obvious for the first individual, and requires a little calculation  
\begin{equation*}
\frac{1}{M-1}\left(1-\frac{1}{M}\right) =  \frac{1}{M}
\end{equation*}
for the second individual. Generally the $k^\text{th}$ individual is introvert with probability
\begin{equation}
\label{Ik}
\frac{1}{M-k+1}\prod_{i=1}^{k-1}\left(1-\frac{1}{M-i+1}\right) = \frac{1}{M}
\end{equation}
Therefore the distribution of introverts is binomial 
\begin{equation}
\label{R:binom}
R(I,M) = \binom{M}{I}\left(1-\frac{1}{M}\right)^{M-I} \left(\frac{1}{M}\right)^I
\end{equation}
from which we conclude that $\langle I\rangle=1$ for any $M$. In contrast to the distributions $P(C,M)$ and $\Pi(F,M)$, the distribution of introverts is asymptotically stationary
\begin{equation}
\label{R:Poisson}
R(I,\infty)  = \frac{e^{-1}}{I!}
\end{equation}

Equation \eqref{F:MRM} yields $n_1(\infty)=\frac{1}{2}$. In Appendix~\ref{ap:RRM}, we show that the degree distribution is asymptotically pure  exponential 
\begin{equation}
\label{nk:RRM}
\lim_{M\to\infty}n_k(M) = 2^{-k}
\end{equation}

When $M$ is large but finite, \eqref{nk:RRM} is applicable as long $\langle N_k(M)\rangle\gg 1$. Alternatively, $k$ should be much smaller than the maximal degree $D(M)$. The maximal degree is found from the criterion \eqref{D:criterion} together with \eqref{nk:RRM} to yield
\begin{equation}
\label{D-max:RRM}
D(M) \simeq \frac{\ln M}{\ln 2}
\end{equation}
When $k\gtrsim \log_2(M)$, fluctuations of $N_k(M)$ become too large and these quantities are no longer self-averaging.

\subsection{Linear choice: Model \eqref{map:linear-rank}}

The average number of egocentric individuals, i.e., the average number of loops is 
\begin{equation}
\label{LM:linear}
\langle \mathcal{L}(M)\rangle = \sum_{i=1}^M \frac{2}{(M-i+1)(M-i+2)}=\frac{2M}{M+1}
\end{equation}

There is again at least one loop as the top rank individual maps to itself, $M\to M$. 
The probability to have one loop 
\begin{equation}
\label{Q1:linear-rank}
Q(1,M) = \prod_{k=2}^{M} \left[1-\frac{2}{k(k+1)}\right] = \frac{M+2}{3M}
\end{equation}
approaches to 1/3 when $M\to\infty$. The probability to have two loops 
\begin{equation}
Q(2,M) = \frac{M+2}{3M} \sum_{j=1}^{M-1} \frac{2}{j(j+3)}
\end{equation}
reduces to
\begin{equation}
\label{Q2:linear-rank}
Q(2,M) = \frac{(M-1)(11M^2+26M+12)}{27M^2(M+1)}
\end{equation}
and approaches to 11/27 when $M\to\infty$. The probability to have three loops is
\begin{equation*}
Q(3,M) = \frac{M+2}{3M} \sum_{1\leq j<k\leq M-1} \frac{4}{j(j+3)k(k+3)}
\end{equation*}
An exact expression of the sum through special functions is rather cumbersome, so we just give
\begin{equation}
\label{Q3:linear-rank}
Q(3,\infty)=  \frac{107 - 6 \pi^2}{243} = 0.196\,635\,282\ldots
\end{equation}
Generalizing above results for small number of loops we arrive at the general formula
\begin{equation*}
Q(\mathcal{L}+1,M) =  \frac{M+2}{3M}  \sum_{1\leq j_1<\ldots<j_\mathcal{L}\leq M-1} \prod_{\nu=1}^\mathcal{L} \frac{2}{j_\nu(j_\nu+3)}
\end{equation*}
Taking the $M\to\infty$ limit and keeping $\mathcal{L}$ finite we obtain 
\begin{equation}
\label{QL:linear-rank}
Q(\mathcal{L}+1,\infty) =  \frac{2^\mathcal{L}}{3} \sum_{1\leq j_1<\ldots<j_\mathcal{L}} \prod_{\nu=1}^\mathcal{L} \frac{1}{j_\nu(j_\nu+3)}
\end{equation}
The probability to have four loops admits a simple expression resembling \eqref{Q3:linear-rank}:
\begin{equation*}
Q(4,\infty)=  \frac{1003 - 90 \pi^2}{2187} = 0.052\,462\,553\ldots
\end{equation*}
We have not found compact formulas for the probabilities $Q(\mathcal{L},\infty)$ with $\mathcal{L}>4$, so we merely mention the asymptotic formula
\begin{equation}
\label{QL:linear-rank-asymp}
Q(\mathcal{L},\infty) \simeq  \frac{2^\mathcal{L}}{(\mathcal{L}-1)!(\mathcal{L}+2)!}
\end{equation}
for $\mathcal{L}\gg 1$.

Returning to finite  $M$ we note that the maximal number of loops occurs with probability
\begin{equation}
\label{QMM:linear}
Q(M,M) =  \frac{2^M}{M!(M+1)!} 
\end{equation}
Equivalently, this is the probability that all individuals are introverts, and the probability $\Pi(0,M)$ that there are no followers. The probability to have the maximal number of followers is 
\begin{equation}
\label{F-max:linear}
\Pi(M-1,M)=\frac{2^M}{(M+1)!} 
\end{equation}

Using the definition of the map, Eq.~\eqref{map:linear-rank}, one finds the average number of followers is given by 
\begin{equation}
\label{Fav:sum}
\langle F\rangle = \langle N_1\rangle  =\sum_{j=1}^M \Psi_j(M)
\end{equation}
where we have used shorthand notation 
\begin{equation}
\label{Psi:def}
\Psi_j(M) = \prod_{i=1}^j \left[1-\frac{2(j-i+1)}{(M-i+1)(M-i+2)}\right]
\end{equation}
Indeed, according to Eq.~\eqref{map:linear-rank}, the individual $j$ is not an expert for $i$ with probability appearing inside the square bracket in the product in \eqref{Psi:def}. The product over all $i\leq j$ ensures that $j$ is not an expert, i.e., a follower.  

An asymptotic analysis of \eqref{Fav:sum}--\eqref{Psi:def} yields 
\begin{equation}
\label{Lav:linear-asymp}
\lim_{M\to\infty}\frac{\langle F\rangle}{M} = n_1 = \frac{e^2-5}{4}= 0.59726402473\ldots
\end{equation}
In Appendix~\ref{ap:derive} we present explicit asymptotic results for $n_k$ with $k\leq 5$. 

Generally, in the $M\to\infty$ limit the degree distribution admits an integral representation
\begin{equation}
\label{nk:linear}
n_{k} = \int_0^1 dx\,\frac{[J(x)]^{k-1}}{(k-1)!}\,e^{-J(x)}
\end{equation}
where
\begin{equation}
\label{J:def}
J(x) = \int_0^x dy\,\frac{2(x-y)}{(1-y)^2} = -2[x + \ln(1-x)]
\end{equation}
An asymptotic analysis of \eqref{nk:linear}--\eqref{J:def} allows one to deduce the large $k$ tail of the degree distribution
\begin{equation}
\label{nk:tail-linear}
n_k\simeq \frac{1}{2 e}\left(\frac{2}{3}\right)^k
\end{equation}
The derivation of \eqref{Lav:linear-asymp}--\eqref{nk:tail-linear} is relegated to Appendix~\ref{ap:derive}.

Using criterion \eqref{D:criterion} together with \eqref{nk:tail-linear} we find that the maximal degree exhibits a logarithmic growth
\begin{equation}
\label{DM:linear}
D(M) \simeq \frac{\ln M}{\ln(3/2)}
\end{equation}

Similarly to \eqref{Ik} we find that the $k^\text{th}$ individual is introvert with probability
\begin{equation}
\label{Ik:linear}
I_k(M)=\frac{2\Psi_k(M)}{(M-k+1)(M-k+2)-2}
\end{equation}
Using the asymptotic behavior of $\Psi_k(M)$ which is derived by employing a continuum approach, see Appendix~\ref{ap:derive}, we deduce
\begin{equation}
\label{Ik:linear-asymp}
I_k(M)\simeq \frac{2}{M^2}\,e^{2k/M}
\end{equation}
The average number of introverts decays as
\begin{equation}
\label{IM:linear} 
\langle I(M)\rangle=\sum_{k=1}^M I_k(M)\simeq  \frac{e^2-1}{M}
\end{equation}

Thus, in the $M\to\infty$ limit, there are on average two egocentrics and no introverts [see \eqref{LM:linear} and \eqref{IM:linear}]. We have also computed the probability distribution of the number of egocentrics, Eq.~\eqref{QL:linear-rank}.

\subsection{Inverse linear choice: Model \eqref{map-linear}}

The average number of egocentric individuals is 
\begin{equation}
\langle \mathcal{L}(M)\rangle = \sum_{k=1}^M \frac{1}{H_k}
\end{equation}
In previous models the average number of egocentric individuals remained finite when $M\to\infty$, the only exception is the model \eqref{map:more} where the number of egocentric individuals diverges, albeit very slowly (logarithmically). In the present model the divergence is much more fast:
\begin{equation}
\label{Lav-linear}
\langle \mathcal{L}(M)\rangle \simeq \frac{M}{\ln M}
\end{equation}
 
The probability to have exactly one egocentric 
\begin{equation}
\label{Q1-linear}
Q(1,M) = \prod_{k=2}^{M} \left[1-\frac{1}{H_k}\right] 
\end{equation}
is very small when $M$ is large
\begin{equation*}
\ln Q(1,M) = -\frac{M}{\ln M} + \mathcal{O}\!\left[\frac{M}{(\ln M)^2}\right]
\end{equation*}
The maximal number of egocentrics occurs with probability
\begin{equation}
\label{QMM:harmonic}
Q(M,M) =  \prod_{k=1}^{M} \frac{1}{H_k}
\end{equation}
which decays faster than exponentially but slower than factorially with $M$:
\begin{equation}
\label{QMM:harm-decay}
Q(M,M) \asymp e^{-M\ln(\ln M)}
\end{equation}
(Here $A \asymp B$ when $M\to\infty$ means the asymptotic equality of the logarithms: $\lim_{M\to\infty} \frac{\ln A}{\ln B}=1$.)

The number of followers is minimal, $F=0$, with probability \eqref{QMM:harmonic}, cf. \eqref{F-min:rank}. The number of followers is maximal, $F=M-1$, with probability 
\begin{equation}
\label{F-max:harmonic}
\Pi(M-1,M) =  \frac{1}{M!}\prod_{j=1}^M\frac{1}{H_j}
\end{equation}

The average number of followers is given by 
\begin{equation}
\label{Fav-linear-exact}
\langle F\rangle = \sum_{j=1}^M \prod_{i=1}^j \left[1-\frac{1}{(j-i+1) H_{M-i+1}}\right]
\end{equation}
The derivation of \eqref{Fav-linear-exact} for the model \eqref{map-linear} is the same as the derivation of \eqref{Fav:sum}--\eqref{Psi:def} for the model \eqref{map:linear-rank}. 

The asymptotic behavior of the average number of followers following from the exact formula \eqref{Fav-linear-exact} is slightly sub-linear, namely
\begin{equation}
\label{Fav-linear}
\langle F\rangle \simeq \frac{M}{\ln M}
\end{equation}
as we show in Appendix~\ref{ap:derive}. Hence in the model \eqref{map-linear} the fraction of followers asymptotically vanishes, although slowly as inverse logarithm: $n_1\simeq (\ln M)^{-1}$. Moreover 
\begin{equation}
\label{nk-linear}
n_k \simeq \frac{1}{\ln M}\,\frac{1}{(k-1)!}
\end{equation}
as we show in Appendix~\ref{ap:derive}. 

As before we find that the $k^\text{th}$ individual is introvert with probability
\begin{equation}
\label{Ik:linear-}
I_k(M)=\frac{1}{H_{M-k+1}}\prod_{j=2}^k \left[1-\frac{1}{j H_{M-k+j}}\right]
\end{equation}
Using this exact formula one finds that the average number of introverts diverges as
\begin{equation}
\label{Iav-linear}
\langle I(M)\rangle \simeq e^{-1}\,\frac{M}{\ln M}
\end{equation}

Thus, in the realm of model \eqref{map-linear}, the number of egocentrics is asymptotically the same as the number of followers [cf. \eqref{Lav-linear} and \eqref{Fav-linear}], and about $37\%$  of egocentrics are introverts [cf. \eqref{Lav-linear} and \eqref{Iav-linear}].

\subsection{The algebraic family of models \eqref{map-a}}

\begin{figure}
\centering
\includegraphics[width=7.77cm]{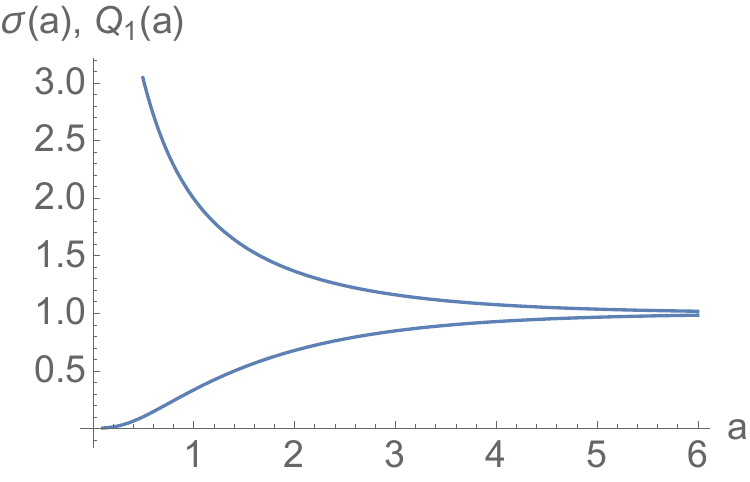}
\caption{Top: The average number of egocentrics $\sigma(a)$ given by \eqref{sigma-a}. Bottom: The probability $Q_1(a)$ to have exactly one egocentric individual given by \eqref{Q1-a}. Shown are results for the one-parameter family of models \eqref{map-a}. }
\label{Fig:sigma}
\end{figure}

The average number of egocentric individuals 
\begin{equation}
\langle \mathcal{L}(M)\rangle = \sum_{k=1}^M \frac{1}{S_k(a)}\,, \qquad S_k(a)=\sum_{p=1}^k p^a
\end{equation}
The average number of egocentric individuals diverges as 
\begin{equation}
\langle \mathcal{L}(M)\rangle \simeq  -\frac{a+1}{aM^a}
\end{equation}
when $-1<a<0$ and saturates to 
\begin{equation}
\label{sigma-a}
\lim_{M\to\infty}\langle \mathcal{L}(M)\rangle = \sigma(a) = \sum_{k=1}^\infty \frac{1}{S_k(a)}
\end{equation}
when $a>0$. The average number of egocentrics $\sigma(a)$ in an infinite random map is a monotonically decaying function of $a$, see Fig.~\ref{Fig:sigma}. Here are a few special values and asymptotic behaviors:
\begin{equation*}
\begin{split}
&\sigma(1) = 2\\
&\sigma(2) = 6(3-4\ln 2) = 1.364\,467\ldots\\
&\sigma(3) = \tfrac{4}{3}\pi^2-12 = 1.159\,472\ldots \\
&\sigma(a) \simeq         a^{-1}    ~~\qquad 0<a\ll 1\\
&\sigma(a) - 1 \simeq 2^{-a}     \quad a\gg 1
\end{split}
\end{equation*}

The probability to have exactly one egocentric is
\begin{equation}
\label{Q1-a}
Q(1,M) = \prod_{k=2}^{M} \left[1-\frac{1}{S_k(a)}\right] 
\end{equation}
The large $M$ behavior of $Q(1,M)$ changes depending on whether $a$ is smaller or bigger than zero:
\begin{equation}
\label{Q1-a-asymp}
Q(1,M) = 
\begin{cases}
\asymp \exp\!\left[-M/\ln M\right]                  & a=-1\\
\asymp \exp\!\left[\frac{a+1}{aM^a}\right]     & -1<a<0\\
M^{-1}                                                           & a=0\\
\simeq Q_1(a)                                               & a>0
\end{cases}
\end{equation}
with 
\begin{equation}
\label{Q1-a-prod}
Q_1(a) = \prod_{k=2}^{\infty} \left[1-\frac{1}{S_k(a)}\right] 
\end{equation}
Here are a few special values  and asymptotic behaviors:
\begin{equation*}
\begin{split}
&Q_1(1) = \frac{1}{3}\\
&Q_1(2) = \frac{3\sqrt{\pi}}{2\Gamma\!\left[\frac{13-\sqrt{-23}}{4}\right] \Gamma\!\left[\frac{13+\sqrt{-23}}{4}\right]} = 0.676\,285\ldots\\
&Q_1(3) = (12\pi)^{-1}\cosh\!\left[\tfrac{\pi\sqrt{7}}{2}\right] = 0.846\,536\ldots \\
&Q_1(a) \asymp e^{-1/a}    \quad\qquad 0<a\ll 1\\
&1-Q_1(a) \simeq 2^{-a} \qquad a\gg 1
\end{split}
\end{equation*}
The probability $Q_1(a)$ to have one egocentric individual is an increasing function of $a$, see Fig.~\ref{Fig:sigma}.

The probability to have two egocentric individuals is 
\begin{equation}
\label{Q2-a}
Q(2,M) = Q(1,M)\sum_{k=2}^{M} \frac{1}{S_k(a)-1}
\end{equation}
When $a>0$, $Q(\mathcal{L},M)\to Q(\mathcal{L},\infty)\equiv Q_\mathcal{L}(a)$ with 
\begin{equation}
\label{QL-a}
Q_\mathcal{L}(a) =  Q_1(a) \sum_{2\leq j_1<\ldots<j_{\mathcal{L}-1}} \prod_{n=1}^{\mathcal{L}-1} \frac{1}{S_{j_n}(a)-1}\\
\end{equation}

The probability that all individuals are egocentrics (equivalently, that there are no followers) is 
\begin{equation}
\label{QMM:Sa}
Q(M,M) =  \prod_{k=1}^{M} \frac{1}{S_k(a)} = \Pi(0,M)
\end{equation}
This formula reduces to \eqref{QMM:harmonic} when $a=-1$, to \eqref{QMM:RRM} when $a=0$, and \eqref{QMM:linear} when $a=1$. The general formula \eqref{QMM:Sa} also simplifies for other integer $a$, e.g., 
\begin{equation*}
Q(M,M) = 
\begin{cases}
\frac{(12)^M}{(M+1)!(2M+1)!}  & a=2\\
\frac{4^M}{[M!(M+1)!]^2}         & a=3
\end{cases}
\end{equation*}

\section{Rewiring}
\label{sec:rewire}

So far, we have analyzed static random maps. Evolving random maps constitute a large and unexplored class of random maps. Let us briefly look at random maps evolving through rewiring: The population size remains fixed while the map varies. For instance, individuals may want to rewire to more connected individuals (i.e., those with higher degrees). We assume that individuals attempting rewiring rely only on local information. Here is a specific example of the rewiring procedure. Starting with a random initial map, we implement rewiring as follows:
\begin{enumerate}
\item At each time step a randomly chosen individual, say $x$, inspects its expert $f(x)$ and the secondary expert $f[f(x)]$. 
\item If the degree of the secondary expert exceeds the degree of the expert, $k(f[f(x)])>k(f(x))$, the switch from $x\to f(x)$ to $x\to f[f(x)]$ occurs; otherwise, the rewiring attempt fails. 
\end{enumerate}
Thus the final map obeys the constraints:
\begin{equation}
k(f[f(x)]) \leq k(f(x)) \quad\text{for all} \quad x\in S
\end{equation}

After a certain number of rewiring steps, the process comes to a halt. For instance, applying the rewiring procedure to the map shown in Fig.~\ref{Fig:random-map-ill} reduces each community to the directed star graph. An example of a directed star graph is shown in Fig.~\ref{Fig:graphs} (bottom). During the rewiring procedure, a follower remains a follower, while an expert may become a follower. Thus the number of followers increases, and the number of experts decreases. If there is a single individual with a maximal degree in a community, the degree of this individual may only increase. Only cycles and directed star graphs (Fig.~\ref{Fig:graphs}) are stable under rewiring. Thus communities of a final map are either cycles and directed star graphs. 

We know the average number of communities in classical random graphs, Eq.~\eqref{c-av}. If the community is a cycle, it remains a cycle; otherwise, it turns into a directed star graph. Each directed star graph has one expert; for a cycle $S_\ell$ of length $\ell$, the number of experts is $\ell$. Thus we have a lower bound $E\geq C$ for the total number of experts after rewiring, leading to 
\begin{equation}
\label{E-av-low}
\langle E\rangle \geq \langle C\rangle = \tfrac{1}{2}\ln M + O(1)
\end{equation}
This lower bound seems qualitatively correct, namely we anticipate a logarithmic growth law 
\begin{equation}
\label{E-av}
\langle E\rangle = B\ln M 
\end{equation}
An obvious challenge is to determine the amplitude $B$.

\section{Experts, Prophets and founders}
\label{sec:prophet}

We call a prophet an expert of the highest degree. For random maps compatible with ranking, the founder is the highest-rank individual. The founder maps to itself, so the founder is an expert. One would like to compute the probability $F_M$ that the founder is a prophet. 

The total number of experts varies within the bounds $1\leq E \leq M$. There is one map with $E=1$, a directed star with a loop from the center to itself; for $M=7$, such map is shown in Fig.~\ref{Fig:graphs} (bottom). Maps with $E=M$ are collections of cycles. The total number of prophets also varies within the bounds 
\begin{equation}
1\leq P \leq M
\end{equation}
Maps with the maximal number of prophets are collections of $M$ loops. The probabilities of such extreme outcome in various random map models are given by \eqref{QMM:classical}, \eqref{QMM:fitness}, \eqref{QMM:0123}, \eqref{QMM:RRM}, \eqref{QMM:linear}, \eqref{QMM:harmonic}, \eqref{QMM:Sa}. 

Little is known about the opposite extreme. The numbers $S(M)$ of maps with a single prophet when all $M^M$ maps are possible have been recently computed via complicated recurrences \cite{Sloane}. The  fraction of maps $S_M/M^M$ with a single prophet is the probability that a classical random map has a single prophet. According to Monte-Carlo simulations, a classical random map has a single prophet with a probability close to $0.35$ for $M=3.2 \cdot 10^6$. In scale-free networks, in contrast, there is a single prophet with probability approaching unity in the $M\to\infty$ limit  \cite{KR13}. 

\begin{figure}
\centering
\includegraphics[width=7.89cm]{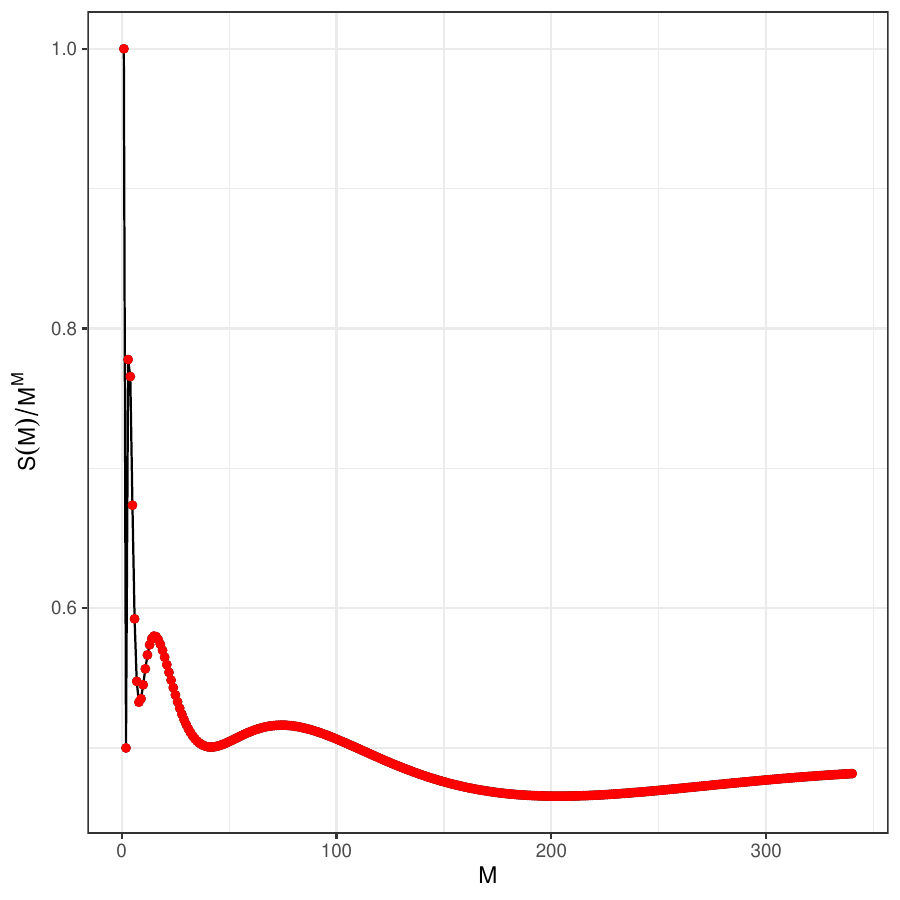}
\caption{The fraction of maps $S_M/M^M$ with a single prophet. Exact results.}
\label{Fig:SM}
\end{figure}

\begin{figure}
\centering
\includegraphics[width=7.89cm]{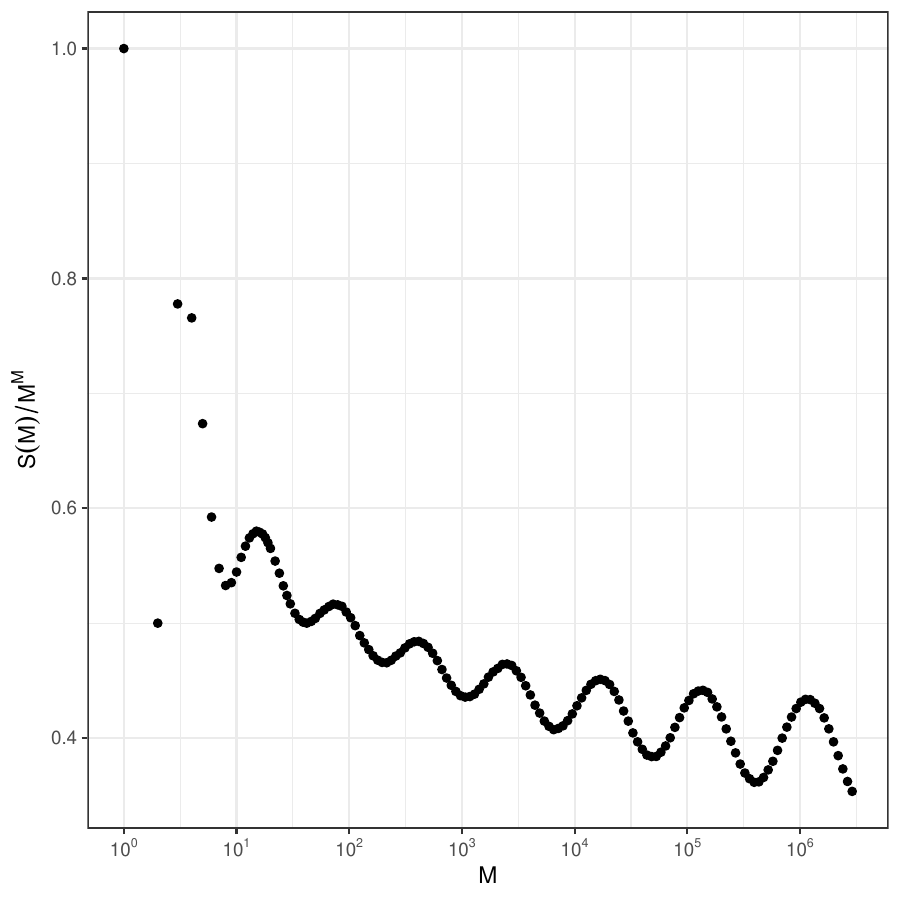}
\caption{The fraction of maps $S_M/M^M$ with a single prophet versus $M$ in a logarithmic scale. Monte-Carlo simulations. }
\label{Fig:SM-log}
\end{figure}

Exact results for $M$ up to a few hundred exhibit intriguing non-monotonic dependence on $M$, see Fig.~\ref{Fig:SM}. Monte-Carlo simulations for $M\leq 3.2 \cdot 10^6$ indicate never-ending oscillations (Fig.~\ref{Fig:SM-log}), with a period of oscillations logarithmically growing with $M$. Heuristic arguments relying on the Poisson form \eqref{nk:Poisson} of the degree distribution of classical random maps suggest a slightly slower growth of period, namely $S(M)/M^M$ is a periodic function of $\ln(M)/\ln[\ln(M)]$. More precisely, the local maxima of $S(M)/M^M$ seemingly  approach a positive constant, while the local minima slowly decay to zero. 

In some random map models, the number of experts is proportional to the population size. The linear scaling $\langle E\rangle\sim M$ is not necessarily a drawback, e.g., if expertise in everyday matters, a more active spouse in a family may be essentially an expert. Still, one anticipates that the number of experts can grow slower than linearly with population size. Perhaps a logarithmic growth is realized for the model with rewiring, cf. \eqref{E-av}. Strictly speaking, we have only established the lower bound \eqref{E-av-low}, so the logarithmic growth \eqref{E-av} has not been proven. 

For classical random maps, the maximal degree is given by \eqref{DM:nu}, i.e., it grows even slightly slower than logarithmically with $M$. Recall that for the family of fitness models \eqref{fitness-a}, the maximal degree grows according to \eqref{DM:nu-a}: if $a\geq 0$, the growth law identical to the growth law \eqref{DM:nu} of the classical random maps; if $-1<a<0$, the growth  is algebraic, $D\sim M^{-a}$; if $a=-1$, a prophet has an almost macroscopic fraction of individuals among the followers, $D\simeq M/\ln M$. 

Characteristics of prophets are sensitive to the type of random map. For models compatible with ranking, the most natural candidate for being a prophet is the founder, i.e., the highest rank individual. For the model \eqref{map:more}, the average degree $\delta$ of the founder and the prophet [cf.  Eq.~\eqref{D-max:RRM}] grow logarithmically with $M$ but with different amplitudes:
\begin{equation}
\label{deg-rank:RRM}
\delta = H_M\simeq \ln M, \quad D(M)\simeq \frac{\ln M}{\ln 2}
\end{equation}
Hence for large $M$, it is unlikely that the prophet is the founder (the first individual in the dynamical interpretation of the RRM). 

For the RRM, one can derive the entire degree distribution $\Phi(k,M)$ of the founder. This distribution satisfies
\begin{equation}
M \Phi(k,M) = \Phi(k-1,M-1) + (M-1) \Phi(k,M-1)
\end{equation}
Solving this recurrence starting with $\Phi(k,1)=\delta_{k,1}$ gives the announced result \eqref{PkM:Stirling}. Using the identity \cite{Knuth}
\begin{equation}
\label{Stirling:GF}
\sum_{k=1}^M {M\brack k} x^k = x(x+1)\cdots(x+M-1)
\end{equation}
one can deduce exact formulas for the average $\delta = \langle k\rangle$ and higher moments. 

In the RRM of the large size $M$, the founder is the prophet with probability
\begin{equation}
\label{Prophet:RRM}
F_M\sim  \sum_{k>D(M)}\Phi(k,M)
\end{equation}
The degree of the prophet is a random variable with small fluctuation when $M\gg 1$. This feature justifies using the average degree $D(M)$ of the prophet in the sum in \eqref{Prophet:RRM}. Using \eqref{PkM:Stirling}, \eqref{Stirling:GF} and the Cauchy theorem, we express $\Phi(k,M)$ as the contour integral
\begin{equation*}
\Phi(k,M) = \frac{1}{2\pi\sqrt{-1}}\oint  \frac{dx}{x^{k}}\,\frac{(x+1)\cdots(x+M-1)}{M!}
\end{equation*}
We  then compute the integral in the $M\to\infty$ limit by applying the saddle point technique. In the interesting situation when $k=\kappa \ln M$ with $\kappa>1$ we get
\begin{equation}
\Phi(k,M) \sim M^{\kappa-1-\kappa\ln\kappa}(\ln M)^{-1/2}
\end{equation}
Combining this with \eqref{Prophet:RRM} and noticing that $k=D(M)$ when $\kappa=1/\ln 2$, we deduce the probability
\begin{equation}
\label{first}
F_M\sim M^{-\epsilon} (\ln M)^{-1/2}, \quad \epsilon = 1-\frac{1+\ln \ln 2}{\ln 2}
\end{equation}
that the founder is the prophet. The exponent $\epsilon$ arises in a surprisingly large number of similarly unrelated problems related to divisors \cite{Divisors,Ford08,Anatomy}, random permutations \cite{Peres16,Ford16a,Ford16b}, averaging processes \cite{BK-21}, etc.  

The RRMs resemble recursive random trees (RRTs) \cite{RRT-Tapia,RRT-Moon,RRT-Pittel}, a null model of growing networks \cite{Newman,book,Frieze} where individuals arrive one by one and map upon arriving to already present individuals. There are no loops in the RRT, but the distinction between the models vanishes as $M\to\infty$, and the chief properties of the maximal degree are expected to be similar. For RRTs, the evolution of the maximal degree was studied in \cite{KR02} where \eqref{first} was derived. Another interesting result of \cite{KR02} is that the expected difference  between the maximal rank $M$ and the rank of the prophet, equivalently, the population size $\mathcal{P}$ when the prophet was born, scales sub-linearly 
\begin{equation}
\label{born}
\mathcal{P} \sim M^\psi, \quad \psi=2-\log_2 3
\end{equation}

For the model \eqref{map:linear-rank}, the average degree of the founder 
\begin{equation}
\delta = \sum_{i=1}^M \frac{2}{M-i+2}= 2[H_{M+1}-1] \simeq 2\ln M
\end{equation}
grows logarithmically with $M$ as the degree of the prophet, Eq.~\eqref{DM:linear}, but the latter grows with larger amplitude $[\ln(3/2)]^{-1}\approx 2.466$. 

For the model \eqref{map-linear}, the average degree of the founder grows much slower than the degree of the prophet 
\begin{equation}
\label{delta-D}
\delta = \sum_{k=1}^M \frac{1}{k H_k} \simeq \ln(\ln M), \quad D(M)\simeq \frac{\ln M}{\ln(\ln M)}
\end{equation}
The degree of the prophet is found from the criterion
\begin{equation*}
1 \sim M\sum_{k>D(M)}n_k\simeq \frac{M/\ln M}{[D(M)]!}
\end{equation*}
where we have used the degree distribution \eqref{nk-linear} and then extracted the asymptotic using the Stirling formula. 

\section{Discussion}
\label{sec:disc}

We proposed random map models with sociological flavor, and explored their basic properties, e.g., the numbers of experts, followers, prophets,  egocentrics, and introverts. Specifically, we analyzed  random maps compatible with fitness or ranking. Some of our results for the average number of egocentrics and introverts in various random map models are collected in Table~\ref{Tab:loops}. The number of egocentrics may increase with $M$. This phenomenon occurs in models \eqref{map:more} and \eqref{map-linear}. In the former model, the number of egocentrics diverges logarithmically with the total number of individuals $M$; in the latter model, it diverges almost linearly, cf. Eq.~\eqref{Lav-linear}. 

In Sect.~\ref{sec:prophet}, we analyzed prophets and founders.  Random maps in which the number of individuals increases provide an attractive setup for studying prophets. One can ask about the size of the population when a prophet was born, the probability that an earlier prophet remains the prophet throughout the evolution, and the total number of distinct prophets that appeared throughout history.  Equations \eqref{first}--\eqref{born} give answers to some of these questions for RRMs. The average number of lead changes increases as $\ln M$ for RRMs, see \cite{KR02}.  The number of distinct prophets is not affected by leap-frogging, so it is smaller than the number of lead changes. The average number of distinct prophets that appeared during the RRM evolution up to size $M$ is unknown.

Random maps compatible with ranking naturally arise in the dynamical setting. Suppose individuals arrive one-by-one, and each new individual $x$ immediately maps to $f(x)$, one of already present individuals (including itself). The RRM model corresponds to the uniform choice. We now define another simple model. Preliminary choice is again uniform, and we postulate that $f(x)=x$ is accepted, but if $f(x)\ne x$, we redirect to the image of the preliminary choice, i.e., the final choice is $f[f(x)]$. We thus obtain recursive random maps with redirection (RRMR). In the RRMR, the community of the founder is a directed star subgraph as in Fig.~\ref{Fig:graphs} (bottom). As the RRM model, the RRMR model is parameter-free and may be tractable. Some characteristics of these two models are identical, e.g., the same formula \eqref{PCM:Stirling} yields the probability distribution of the number of communities. Other characteristics differ and are sometimes easier to determine for the RRMR. For instance, the degree distribution $\Phi(k,M)$ of the founder of the RRMR satisfies
\begin{equation*}
M \Phi(k,M) = (k-1)\Phi(k-1,M-1) + (M-k) \Phi(k,M-1)
\end{equation*}
Solving this recurrence starting with $\Phi(k,1)=\delta_{k,1}$ we find that the degree distribution of the founder is flat
\begin{equation}
\label{Phi:RRMR}
\Phi(k,M) = M^{-1}, \qquad k = 1, \ldots, M
\end{equation}
Using \eqref{Phi:RRMR}, we find the average degree of the founder for the RRMR, $\delta = (M+1)/2$, significantly exceeds the average degree of the founder  for the RRM, $\delta=H_M$.

The dynamical setting may help to appropriate the techniques developed for studying similar types of growing random trees \cite{KR02,JML08,JML10} to explore more subtle questions about prophets and the founder for the RRMR and other dynamical random maps compatible with ranking.

Vertices of anomalously high degrees are widespread in modern communication and social networks. Referring to such vertices as prophets is questionable. Intuitively, a prophet is more like a central vertex in the largest community, and the size of the community appears a more interesting quantity than the degree of the prophet. Thus, rather than studying the degrees of prophets and founder (in random maps compatible with ranking), one could investigate the largest community and the primordial community (containing the founder). 

The size distribution of the largest community is difficult to  probe analytically. In most random map models, the distribution $\Psi(s,M)$ of the size $s$ of the largest community acquires the scaling form
\begin{equation}
\Psi(s,M) = M^{-1}\psi(\sigma)
\end{equation}
when $s\to\infty$ and $M\to\infty$ with $\sigma = s/M$ kept finite. The distribution $\psi(\sigma)$ has singularities at $\sigma=1/m$ with $m=2,3,\ldots$ that become weaker and weaker as $m$ increases. The origin of these singularities is easy to appreciate: One community may have size exceeding $M/2$, i.e., in the $\frac{1}{2}<\sigma<1$ range; two communities may have size in the $\frac{1}{3}<\sigma<\frac{1}{2}$ range; etc. This infinite set of singularities prevents an analytical determination of $\psi(\sigma)$. 
 
For classical random maps, the rescaled size distribution $\psi(\sigma)$ of the largest community was probed numerically in \cite{BN-Derrida87}. The average size of the largest community is close to $\lambda M$ with (see \cite{Flajolet90})
\begin{equation}
\label{lambda:RM}
\lambda = 2\int_0^\infty dy\,\left[1-e^{-\Gamma(0,y)/2}\right]=0.757\,823\ldots
\end{equation}

For the RRM, the rescaled size distribution and the average size of the largest community are unknown. The primordial community keeps its identity throughout the evolution. The size distribution $\mathcal{P}(s,M)$ of the primordial community satisfies 
\begin{equation*}
\mathcal{P}(s,M+1) = \tfrac{s-1}{M+1}\,\mathcal{P}(s-1,M) +\tfrac{M+1-s}{M+1}\,\mathcal{P}(s,M)
\end{equation*}
from which we find that the size distribution of the primordial community is flat
\begin{equation}
\label{P:RRM}
\mathcal{P}(s,M) = M^{-1}, \qquad s = 1, \ldots, M
\end{equation}
Therefore, the average size of the primordial community is $\langle s \rangle = \frac{M+1}{2}$. 

The average size of the largest community exceeds the average size of the primordial community and also scales linearly with $M$, namely, as $\lambda M$ with $\lambda>\frac{1}{2}$. Computing the amplitude $\lambda$ is perhaps feasible as it was  \cite{Flajolet90} for classical random maps, Eq.~\eqref{lambda:RM}.

\bigskip\noindent
I benefitted from discussions with Mirta Galesic and Sid Redner. I am particularly grateful to Aaron Schweiger, the author of the sequence A365061 in \cite{Sloane}, for supplying Figs.~\ref{Fig:SM}--\ref{Fig:SM-log}, and for helpful comments.

\appendix

\section{Recursive random maps (RRMs)}
\label{ap:RRM}

Here, we consider an evolving random map with a growing number of individuals. We postulate that individuals arrive one by one, and each new individual $x$ immediately maps to an individual that we denote $f(x)$. As in classical random maps, we postulate that the image $f(x)$ is chosen uniformly and randomly among all present individuals, including $x$. The outcome is a recursive random map (RRM) model. 

If the choice is uniform but excludes the new individual, $f(x)\ne x$, we almost recover recursive random trees \cite{RRT-Tapia,RRT-Moon,RRT-Pittel}, a null model of growing networks \cite{Newman,book,Frieze}. The only exception is the founder, which must map to itself. These RRTs with an extra loop contain a single community and one egocentric, the founder, so they are less natural in the sociological realm than the RRM. A model of growing random hypergraphs known as random recursive hypergraphs \cite{PK-hyper} also overlaps with RRMs.

For the RRM, the number of communities $C(M)$ is a random variable that grows according to
\begin{equation}
\label{CM}
C(M) = 
\begin{cases}
C(M-1)               & \text{prob} ~~ 1-M^{-1}\\
C(M-1) + 1         & \text{prob} ~~ M^{-1}
\end{cases}
\end{equation}
In particular, the average satisfies the recurrence
\begin{equation}
\label{CM-av}
\langle C(M) \rangle = \langle C(M-1) \rangle + M^{-1}
\end{equation}
which is solved to yield
\begin{equation}
\label{CM-av-sol}
\langle C(M) \rangle = H_M
\end{equation}
Comparing with \eqref{c-av}, we see that there are, on average, twice as many communities in the RRM than in classical random maps.

More generally, using \eqref{CM} we find that the probability distribution defined by Eq.~\eqref{PCM:def} satisfies
\begin{eqnarray}
\label{PCM}
P(C,M) &=& (1-M^{-1})P(C,M-1) \nonumber \\
             &+& M^{-1}P(C-1,M-1)
\end{eqnarray}
The initial condition is
\begin{equation}
\label{P:IC}
P(C,1)=\delta_{C,1}
\end{equation}

When $C=1$, Eq.~\eqref{PCM} reduces to
\begin{equation*}
P(1,M) = (1-M^{-1})P(1,M-1)
\end{equation*}
Solving this recurrence subject to $P(1,1)=1$ yields
\begin{equation}
\label{P1M-sol}
P(1,M) = M^{-1}
\end{equation}
Using \eqref{P1M-sol} we recast \eqref{PCM} for $C=2$ into
\begin{equation*}
\label{P2M}
P(2,M) = (1-M^{-1})P(2,M-1)+\frac{1}{M(M-1)}
\end{equation*}
Solving this recurrence subject to $P(2,1)=0$ yields
\begin{equation}
\label{P2M-sol}
P(2,M) = M^{-1}H_{M-1}
\end{equation}
for $M\geq 2$.  Collecting \eqref{P1M-sol} and \eqref{P2M-sol} give \eqref{P12}. 

We now derive the announced general solution \eqref{PCM:Stirling}. We multiply \eqref{PCM} by $M!$ and notice that $M! P(C,M)$ satisfies an addition formula for the Stirling numbers ${M\brack C}$ of the first kind \cite{Knuth}. The initial condition \eqref{P:IC} agrees with the standard initial condition for the Stirling numbers of the first kind, and we arrive at \eqref{PCM:Stirling}. 

Using \eqref{PCM:Stirling} and identities for the relevant Stirling numbers of the first kind, we extract the probabilities that the number of communities is close to maximal:
\begin{equation}
\label{PMM3}
\begin{split}
 P(M,M)  &= \frac{1}{M!}  \\
P(M-1,M) & = \frac{1}{M!}\binom{M}{2} \\
 P(M-2,M)  &= \frac{1}{M!}\,\frac{3M-1}{4}\binom{M}{3} \\
 P(M-3,M)  &= \frac{1}{M!}\binom{M}{2}\binom{M}{4}
\end{split}
\end{equation}
etc. In the complementary limit of $C=O(1)$ and $M\gg 1$
\begin{equation}
\label{PCM:log}
P(C,M)  \simeq  M^{-1}\,\frac{(\ln M)^{C-1}}{(C-1)!}
\end{equation}

In the RRM model, every new community starts as a loop, so the total number of communities is equal to the total number of loops: $C=\mathcal{L}$. This is not applicable to all maps, e.g., for the map shown in Fig.~\ref{Fig:random-map-ill} we have $C=3$ and $\mathcal{L}=1$. For RRMs, however, the equality $C=\mathcal{L}$ always holds. Thus \eqref{P12} imply
\begin{equation*}
Q(1,M) =  M^{-1}, \quad Q(2,M) = M^{-1}H_{M-1}
\end{equation*}
and one can similarly re-express \eqref{CM-av-sol} and \eqref{PMM3}--\eqref{PCM:log}. 

These results are identical to corresponding results for {\em static} random maps compatible with ranking and satisfying \eqref{map:more}. To establish an isomorphism between the dynamic RRM model and the static model \eqref{map:more}, we label individuals in the dynamic RRM model in reverse order. Indeed, if we grow the random map to size $M$, we assign rank $M$ to the first individual, rank $M-1$ to the second individual, etc.; the equivalence of the rules is then easy to appreciate.

The number of followers $F(M)$ is a random variable evolving according to 
\begin{equation}
\label{FM:eq}
F(M+1) =
\begin{cases}
F(M)       & \text{prob} ~~\frac{F(M)+1}{M+1}\\
F(M) + 1 & \text{prob} ~~\frac{M-F(M)}{M+1}
\end{cases}
\end{equation}
Indeed, the number of followers does not change if the newly-introduced individual maps to a follower or itself. This happens with probability $\frac{F(M)+1}{M+1}$. Otherwise, the number of followers increases by one. The stochastic rule \eqref{FM:eq} implies that the probability distribution defined by Eq.~\eqref{PFM:def} satisfies
\begin{eqnarray}
\label{PFM}
M\Pi(F,M) &=& (F+1)\Pi(F,M-1) \nonumber \\
             &+& (M-F)\Pi(F-1,M-1)
\end{eqnarray}

Multiplying Eq.~\eqref{PFM} by $(M-1)!$ and noticing that $M! \Pi(F,M)$ satisfies recurrence identical to the addition formula for Eulerian numbers \cite{Euler36,Euler55,Knuth} we arrive at the announced result \eqref{PFM:Euler}. Explicit expressions \cite{Knuth}  
\begin{equation*}
\begin{split}
\eulerian{M}{0} &=1 \\
\eulerian{M}{1}  &=2^M-M-1 \\
\eulerian{M}{2} &=3^M-(M+1)2^M+\binom{M+1}{2} \\
\eulerian{M}{3} &=4^M-(M+1)3^M+2^M\binom{M+1}{2}-\binom{M+1}{3}
\end{split}
\end{equation*}
for Eulerian numbers give $\Pi(F,M)$ for $F=0,1,2,3$. While $\Pi(0,M)=1/M!$ is obvious from the definition of the RRM, a straightforward derivation of $\Pi(F,M)$ for $F=1,2,3$ is laborious. Explicit results for $\Pi(F,M)$ with $F=M-1,M-2,M-3, M-4$ then follow from the mirror symmetry between Eulerian numbers \cite{Knuth}
\begin{equation}
\label{mirror:E}
\eulerian{M}{F}=\eulerian{M}{M-1-F}
\end{equation}

We now derive the announced asymptotic behavior \eqref{nk:RRM} of the degree distribution.  The number of vertices of
degree $k\geq 2$ evolves according to stochastic rules
\begin{equation*} 
\label{NkM}                                                                               
N_k(M+1)=
\begin{cases}
N_k(M)-1    & \mathrm{prob}~~ \frac{N_k(M)}{M+1}\\[2mm]
N_k(M)+1   & \mathrm{prob}~~ \frac{N_{k-1}(M)}{M+1}\\[2mm]
N_k(M)       & \mathrm{prob}~~ 1-\frac{N_{k-1}(M)+N_k(M)}{M+1}
\end{cases}                                       
\end{equation*}
Averaging above equation yields the recurrence
\begin{equation}
\label{Nkav}
(M+1)\langle N_k(M+1)\rangle
=M\langle N_k(M)\rangle +\langle N_{k-1}(M)\rangle
\end{equation}
When $k=2$, this recurrence gives
\begin{equation}
\label{N2av}
(M+1)\langle N_2(M+1)\rangle
=M\langle N_2(M)\rangle + \frac{M-1}{2}
\end{equation}
where we have used 
\begin{equation*}
\left\langle N_1(M)\right\rangle = \langle F(M)\rangle = \frac{M-1}{2}
\end{equation*}
The initial graph is a loop. Solving \eqref{N2av} subject to the initial condition $N_2(1)=1$ yields
\begin{equation}
\label{n2M}
n_2(M) = \frac{1}{4}-\frac{3}{4}M^{-1}+\frac{3}{2}M^{-2}
\end{equation}

Specializing  \eqref{Nkav} to $k=3$ and using \eqref{n2M} we obtain
\begin{equation*}
\label{N3av}
(M+1)\langle N_3(M+1)\rangle = M\langle N_3(M)\rangle + \frac{M}{4} -\frac{3}{4} +\frac{3}{2}M^{-1}
\end{equation*}
which we solve subject to $N_3(1)=0$ to yield
\begin{equation}
\label{n3M}
n_3(M) = \frac{1}{8}-\frac{7}{8}M^{-1}+\frac{3+6H_{M-1}}{4M^2}
\end{equation}

Exact results for $n_k(M)$  become more and more cumbersome as $k$  increases. The most important leading and sub-leading terms are simple. In the $M\to\infty$ limit, the recurrence \eqref{Nkav} simplifies to $2n_k=n_{k-1}$ from which we deduce \eqref{nk:RRM}. A more careful asymptotic analysis of \eqref{Nkav} allows one to extract the sub-leading term:
\begin{equation}
\label{nkM}
n_k(M) = 2^{-k}-\big(1-2^{-k}\big)M^{-1}+\ldots
\end{equation}

\section{Derivation of \eqref{Lav:linear-asymp}--\eqref{nk:tail-linear} and \eqref{Fav-linear}--\eqref{nk-linear}}
\label{ap:derive}

We begin with model \eqref{map:linear-rank}. To compute the sum in Eq.~\eqref{Fav:sum} we first simplify $\Psi_j(M)$ defined by \eqref{Psi:def}. Taking the logarithm of the product in \eqref{Psi:def} we obtain 
\begin{equation*}
\ln \Psi_j(M) = \sum_{k=1}^j \ln\left[1-\frac{2k}{(M-j+k)(M-j+k+1)}\right]
\end{equation*}
where $k=j-i+1$. Replacing summation on the right-hand side by integration (this is justifiable since we are interested in the $M\to\infty$ behavior) we arrive at
\begin{equation*}
\ln \Psi_j(M) \simeq \int_0^j dk\,\ln\!\left[1-\frac{2k}{(M-j+k)^2}\right]
\end{equation*}
in the leading order. Expanding the logarithm and keeping the leading term we obtain
\begin{equation*}
\frac{\ln \Psi_j(M)}{2} \simeq  -\int_0^j \frac{k\,dk}{(M-j+k)^2} = \frac{j}{M} + \ln\!\left(1-\frac{j}{M}\right)
\end{equation*}
Plugging this into \eqref{Fav:sum} and replacing  summation over $j$ by integration over  $x=j/M$ 
we obtain
\begin{equation}
\label{Lav:int}
\langle F\rangle \simeq  M\int_0^1 dx\, (1-x)^2 e^{2x}=M\int_0^1 dx\,e^{-J(x)}
\end{equation}
with $J(x)$ given by \eqref{J:def}.  

Similarly to $\langle F\rangle = \langle N_1\rangle$ one finds
\begin{equation*}
\langle N_2\rangle = \sum_{j=1}^M \sum_{i=1}^j  \frac{2(j-i+1)\,\Psi_j(M)}{(M-i+1)(M-i+2)-2(j-i+1)}
\end{equation*}
Replacing summations by integrations one gets
\begin{equation}
\label{n2:linear}
n_2 = \int_0^1 dx\, J(x)\,e^{-J(x)}
\end{equation}

Similarly one derives a general formula for $\langle N_k\rangle$ involving $k$ sums. Summations can be replaced by integrations when $M\to\infty$. This leads to the announced Eq.~\eqref{nk:linear}. 

One can explicitly compute $n_k$ for small $k$. The integral in \eqref{Lav:int} is $(e^2-5)/4$ leading to the announced result \eqref{Lav:linear-asymp}. Computing the integral in \eqref{n2:linear} we obtain
\begin{equation}
\label{n2:linear-sol}
n_2 = -1+ \tfrac{1}{2}e^2 C, \quad C = \gamma + \ln 2 - 1 + \int_{2}^\infty \frac{dt}{t}\,e^{-t}
\end{equation}
where we have introduced the quantity  $C$ that also appears in $n_k$ for larger $k$. Indeed, computing the integrals in Eq.~\eqref{nk:linear} one gets
\begin{equation}
\label{n345}
\begin{split}
&n_3 = 2e^2S_3 -\tfrac{5e^2+1}{4}-e^2 C\\
&n_4 = 4e^2(S_4-S_3) - \tfrac{e^2+9}{12} -  \tfrac{1}{2}\,e^2 C  \\
&n_5 = 2e^2(4S_5 - 4S_4 - S_3) + \tfrac{13e^2-11}{12} +  \tfrac{5}{6}\,e^2 C
\end{split}
\end{equation}
etc. where we shortly write 
\begin{equation}
S_p \equiv \sum_{j\geq 0} \frac{1}{(j+1)^p}\,\frac{(-2)^j}{j!}
\end{equation}
Here are a few numerical values:
\begin{equation*}
\begin{split}
&n_1 = 0.59726402473266\ldots \\
&n_2 = 0.17952742450736\ldots \\
&n_3 = 0.08878616080733\ldots \\
&n_4 = 0.05031763625502\ldots \\
&n_5 = 0.03034754225441\ldots 
\end{split}
\end{equation*}

To derive Eq.~\eqref{nk:tail-linear} we first replace the integration variable in Eq.~\eqref{nk:linear} from $x$ to $J$:
\begin{equation}
\label{nk:linear-J}
n_{k} = \int_0^\infty dJ\,\frac{J^{k-1}}{(k-1)!}\,e^{-J}\,\frac{1-x}{2x}
\end{equation}
When $k\gg 1$, the dominant contribution to the integral comes from large $J$, so Eq.~\eqref{J:def} gives $1-x\simeq e^{-1-J/2}$, and \eqref{nk:linear-J} simplifies to the announced asymptotic \eqref{nk:tail-linear}:
\begin{equation}
\label{nk:linear-asymp}
n_{k} \simeq \frac{1}{2e}\int_0^\infty dJ\,\frac{J^{k-1}}{(k-1)!}\,e^{-3J/2}=\frac{1}{2e}\left(\frac{2}{3}\right)^k
\end{equation}

We now turn to model \eqref{map-linear}. We use \eqref{Fav:sum} with
\begin{equation}
\label{Psi:2}
\Psi_j(M) = \prod_{i=1}^j \left[1-\frac{1}{(j-i+1)H_{M-i+1}}\right]
\end{equation}
Taking the logarithm, we obtain the sum
\begin{equation}
\label{Pj-M}
\ln \Psi_j(M) = \sum_{k=1}^j \ln\left[1-\frac{1}{kH_{M-j+k}}\right]
\end{equation}
where $k=j-i+1$. The asymptotic behavior of \eqref{Pj-M} is
\begin{equation}
\label{P-HH}
\ln \Psi_j(M) \simeq -\frac{H_j}{H_M}
\end{equation}
Combining \eqref{Fav:sum} and \eqref{P-HH} we obtain
\begin{eqnarray*}
\langle N_1\rangle &=&  \sum_{j=1}^M  \Psi_j(M) \\
& \simeq &  \int_1^M dj\,\exp\!\left[-\frac{\ln j}{\ln M}\right] 
\simeq  \int_1^M dj\,\left[1-\frac{\ln j}{\ln M}\right] 
\end{eqnarray*}
leading to \eqref{Fav-linear}. Similarly we find
\begin{eqnarray*}
\langle N_2\rangle &=& \sum_{j=1}^M \sum_{i=1}^j  \frac{\Psi_j(M)}{(j-i+1)H_{M-i+1}-1} \\
&\simeq& \sum_{j=1}^M \frac{\Psi_j(M) H_j}{H_M} \simeq \int_1^M dj\,\left[1-\frac{\ln j}{\ln M}\right] \frac{\ln j}{\ln M}
\end{eqnarray*}
leading to $n_2\simeq (\ln M)^{-1}$. 

A general formula for $\langle N_k\rangle$ involves $k$ sums. Simplifying it as above and replacing summations by integrations one establishes \eqref{nk-linear}.

\bibliography{references-nets}

\end{document}